\documentclass[11pt]{article}

\usepackage{amsmath}
\usepackage{amsmath, amsthm, amsfonts, amssymb, color}
\usepackage{mathrsfs}
\usepackage{latexsym,bm}
\usepackage{graphicx}

\newtheorem{theorem}{Theorem}[section]%edit the theorems
\newtheorem{proposition}[theorem]{Proposition}

\newtheorem{lemma}[theorem]{Lemma}

\newtheorem{remark}[theorem]{Remark}

\newcommand{\me}{\mathcal{M}(\mathbb{R})}
\newcommand{\p}{\mathrm{P}}% for skeleton space
\newcommand{\pp}{\mathbb{P}}% for the original super-BM
\newcommand{\bp}{\mathbf{P}}% for the skeleton branching diffusion
\newcommand{\Bp}{\Pi}% for the Brownian motion
\newcommand{\mc}{\mathcal{M}_{c}(\mathbb{R})}
\newcommand{\mf}{\mathcal{M}(\mathbb{R})}

\newcommand{\B}{\mathcal{B}}
\newcommand{\R}{\mathbb{R}}

\newcommand{\1}{\mathbf{1}}
\newcommand{\e}{\mathrm{e}}

\definecolor{wco}{rgb}{0.5,0.2,0.3}

\numberwithin{equation}{section} %the number of the equation will be counted in one section
\begin{document}

\allowdisplaybreaks

\title{\bf Limiting distributions for a class of super-Brownian motions with spatially dependent branching mechanisms
\footnote{The research of this project is supported
     by the National Key R\&D Program of China (No. 2020YFA0712900).}}
\author{ \bf  Yan-Xia Ren\footnote{The research of this author is supported by  NSFC (Grant Nos. 12071011 and 12231002) and  the Fundamental Research Funds for Central Universities, Peking University LMEQF.\hspace{1mm} } \hspace{1mm}\hspace{1mm}
Ting Yang\footnote{ The research of this author is supported by NSFC (Grant No. 12271374).}
\hspace{1mm} }
\date{}
%\date{}
\maketitle

\begin{abstract}
In this paper we consider a large class of super-Brownian motions in $\R$ with spatially dependent branching mechanisms. We establish the almost sure growth rate of the
mass located outside a time-dependent interval $(-\delta t,\delta t)$ for $\delta>0$. The growth rate is given in terms of the principal eigenvalue $\lambda_{1}$ of the Sch\"{o}dinger type operator associated with the branching mechanism. From this result we see the existence of phase transition for the growth order at $\delta=\sqrt{\lambda_{1}/2}$. We further show that the super-Brownian motion shifted by
$\sqrt{\lambda_{1}/2}\,t$
converges in distribution to a random measure with random density mixed by a martingale limit.
\end{abstract}

\bigskip


\bigskip


\section{Introduction and main results}\label{sec:prelim}

\subsection{Super-Brownian motions}

Let $\mathcal{M}(\R)$ (resp. $\mathcal{M}_{c}(\R)$) denote the set of finite (resp. finite and compactly supported) measures on $\R$. When $\mu$ is a measure on $\R$ and $f$ is a measurable function, define $\langle f,\mu\rangle=\int_{\R}f(x)\mu(dx)$ whenever the right hand side makes sense. Sometimes we also write $\mu(f)$ for $\langle f,\mu\rangle$.
Let $\left((B_{t})_{t\ge 0},\Bp_{x},x\in\R \right)$
be a standard Brownian motion on $\R $ with $\Bp_{x}\left(B_{0}=x\right)=1$.
The main process of interest in this paper is an $\me$-valued Markov process $X=\{X_{t}:t\ge 0\}$ with evolution depending on two quantities $P_{t}$ and $\psi$. Here $P_{t}$ is the semigroup of
$\left((B_{t})_{t\ge 0},\Bp_{x},x\in\R \right)$
and $\psi$ is the so-called
branching mechanism,
 which takes the form
\begin{equation}\label{bm}
\psi(x,\lambda)=-\beta(x)\lambda+\alpha(x)\lambda^{2}+\int_{(0,+\infty)}\left(\e^{-\lambda u}-1+\lambda u\right)\pi(x,du)\quad x\in\R ,\lambda\ge 0,
\end{equation}
where $\beta\in C_{c}(\R )$, $0\not\equiv\alpha\in C^{+}_{c}(\R )$, and $\pi$ is a kernel from $\R$ to $(0,+\infty)$ such that
$$\int_{(0,+\infty)}u^{2}\pi(x,du)\in C^{+}_{c}(\R).$$
The distribution of $X$ is denoted by $\pp_{\mu}$ if it is started
at $\mu\in\me$ at $t=0$.
$X$ is called a $(B_{t},\psi)$-superprocess or super-Brownian motion with branching mechanism $\psi$ if for all $\mu\in \me$, nonnegative bounded measurable function $f$ and $t\ge 0$,
 \begin{equation}
 \pp_{\mu}\left[e^{-\langle f,X_{t}\rangle}\right]=e^{-\langle u_{f}(t, \cdot),\mu\rangle},\label{eq0}
 \end{equation}
 where $u_{f}(t,x)=-\log \pp_{\delta_{x}}\left(e^{-\langle f,X_{t}\rangle}\right)$ is
 the unique nonnegative locally bounded solution
 to the following integral equation:
 \begin{equation}\label{eq1}
 u_{f}(t,x)=P_{t}f(x)-\int_{0}^{t}P_{s}\left(\psi(\cdot,u_{f}(t-s, \cdot))\right)(x)ds\quad\forall x\in\R ,\ t\ge 0.
 \end{equation}
The existence of such a process $X$ is established in \cite{D93}.
A closely related  $\me$-valued  process is branching Brownian motion
with the branching rate given by either a compactly supported measure or a function decaying sufficiently fast at infinity (see, e.g., \cite{Bocharov,BH,Shiozawa1,Shiozawa2,NS,Nishimori} and references therein).

\subsection{Notation and some facts}

Throughout this paper we use ``$:=$" to denote a definition.
For functions $f$ and $g$ on $\R $, $\|f\|_{\infty}:=\sup_{x\in \R }|f(x)|$ and $(f,g):=\int_{-\infty}^{+\infty }f(x)g(x)dx$.
For positive functions $f(x)$ and $g(x)$ on $(0,+\infty)$, we write $f(x)\sim g(x) \, (x\to+\infty)$ if $\lim_{x\to+\infty}f(x)/g(x)=1$.
For $a,b\in \mathbb{R}$, $a\wedge b:=\min\{a,b\}$, $a\vee b:=\max\{a,b\}$.
The letters $c$ and $C$ (with subscript) denote finite positive constants which may vary from place to place.

Let $\mathcal{M}_{loc}(\R)$ denote the space of locally finite Borel measures on $\R$
with vague topology,
which is generated by the integration maps
 $\pi_{f}:\mu\mapsto\mu(f)$
 for all compactly supported bounded continuous functions $f$ on $\R$. A random variable taking values in $\mathcal{M}_{loc}(\R)$ is called a random measure on $\R$. We say random measures $\xi_{n}$ converges in distribution to $\xi$ if $\mathrm{E}[F(\xi_{n})]\to \mathrm{E}[F(\xi)]$ for every bounded continuous function $F$ on $\mathcal{M}_{loc}(\R)$. \cite[Theorem 4.11]{Kallenberg} proves that $\xi_{n}$ converges in distribution to $\xi$ if and only if the random variables $\langle f,\xi_{n}\rangle$ converges in distribution to $\langle f,\xi\rangle$ for every $f\in C^{+}_{c}(\R)$.

For a measurable function $f$, we set
$$e_{f}(t):=\exp\left\{\int_{0}^{t}f(\xi_{s})ds\right\},\quad t\ge 0,$$
whenever it is well defined.
We define the Feynman-Kac semigroup $P^{\beta}_{t}$ by
 \begin{equation}
 P^{\beta}_{t}f(x):=\Bp_{x}\left[e_{\beta}(t)
 f(\xi_{t})\right]\quad\mbox{ for  }
 f\in\mathcal{B}^{+}_{b}(\R).\nonumber
 \end{equation}
 Define
  \begin{equation}\label{def-gamma}
  \gamma(x):=\alpha(x)+\frac{1}{2}\int_{(0,+\infty)}u^{2}\pi(x,du),\quad x\in \R.
  \end{equation}
It is known (cf. \cite{D93}) that for every $\mu\in\mf$ and $f\in\mathcal{B}^{+}_{b}(\R)$, the first two moments of $\langle f,X_{t}\rangle$ exist and can be expressed as
 \begin{equation}\label{mean}
 \pp_{\mu}\left(\langle f,X_{t}\rangle\right)=\langle P^{\beta}_{t}f,\mu\rangle,
 \end{equation}
and
\begin{equation}\label{var}
\mbox{Var}_{\mu}\left(\langle f,X_{t}\rangle\right)=\int_{0}^{t}\langle P^{\beta}_{s}\left(2\gamma\left(P^{\beta}_{t-s}f\right)^{2}\right),\mu\rangle ds.
\end{equation}

The spectrum of the operator $\mathcal{L}=\frac{1}{2}\Delta+\beta$, denoted by $\sigma(\mathcal{L})$, consists of $(-\infty,0]$ and at most a finite number of nonnegative eigenvalues.
Throughout this paper, we make the following assumption:
\begin{equation}\label{a1}
\lambda_{1}:=\sup(\sigma(\mathcal{L}))>0.\tag{A1}
\end{equation}
Then $\lambda_{1}$ is simple and the corresponding eigenfunction (ground state) $h$ can be taken to be strictly positive, bounded and continuous. We choose $h$ that is normalized with $\int_{-\infty}^{+\infty}h^{2}(x)dx=1$.
We remark here that \eqref{a1} is automatically satisfied when $\beta\ge 0$ is a nontrivial function.

One has (see, for example, \cite[Lemma 3.1]{NS})
\begin{equation}\label{eq:h}
h(x)=\int_{-\infty}^{+\infty}G_{\lambda_{1}}(x,y)\beta(y)h(y)dy.
\end{equation}
where $G_{\lambda_{1}}(x,y)$ denotes
the $\lambda_{1}$-potential density
of Brownian motion. Using the fact that
\begin{equation}\label{estiforG}
G_{\lambda_{1}}(x,y)\sim \frac{1}{\sqrt{2\lambda_{1}}}\mathrm{e}^{-\sqrt{2\lambda_{1}}|x-y|}\mbox{ as }|x-y|\to +\infty,
\end{equation}
one can easily show that
\begin{equation}\label{estimateforh}
h(x)\sim C_{\mp}\mathrm{e}^{-\sqrt{2\lambda_{1}}|x|}\mbox{ as }x\to \pm\infty,
\end{equation}
where
\begin{equation}\label{def-C}C_{\mp}:=\frac{1}{\sqrt{2\lambda_{1}}}\int_{-\infty}^{+\infty}\beta(y)h(y)\mathrm{e}^{\pm\sqrt{2\lambda_{1}}y}dy.\end{equation}

Since $\e^{-\lambda_{1}t}P^{\beta}_{t}h=h$ for all $t\ge 0$, one can show by the Markov property that
$$W^{h}_{t}(X):=\e^{-\lambda_{1}t}\langle h,X_{t}\rangle,\quad \forall t\ge 0$$
is a nonnegative $\pp_{\mu}$-martingale for every $\mu\in\me$. Let $W^{h}_{\infty}(X)$ be the martingale limit.
It then follows by \cite[Theorem 3.2]{PY} that for every nontrivial $\mu\in\mc$,
$$\lim_{t\to+\infty}W^{h}_{t}(X)=W^{h}_{\infty}(X)\quad \pp_{\mu}\mbox{-a.s. and in }L^{2}(\pp_{\mu}).$$
Hence $W^{h}_{\infty}(X)$ is non-degenerate in the sense that $\pp_{\mu}\left(W^{h}_{\infty}(X)>0\right)>0$.

\subsection{Main results}

For any $R\ge 0$, define
$$\mathcal{X}^{R}_{t}:=\langle \1_{(-R,R)^{c}},X_{t}\rangle.$$

\begin{theorem}\label{them1}
For any $\delta>\sqrt{\lambda_{1}/2}$ and $\mu\in\mc$,
$$\lim_{t\to+\infty}\mathcal{X}^{\delta t}_{t}=0\quad\pp_{\mu}\mbox{-a.s.}$$
For any $0\le \delta<\sqrt{\lambda_{1}/2}$ and $\mu\in\mc$,
$$\lim_{t\to+\infty}\frac{\log \mathcal{X}^{\delta t}_{t}}{t}=\lambda_{1}-\sqrt{2\lambda_{1}}\delta\quad \pp_{\mu}\mbox{-a.s. on }\{W^{h}_{\infty}(X)>0\}.$$
\end{theorem}

According to Theorem \ref{them1},
for $\delta<\sqrt{\lambda_{1}/2}$,
the mass outside $(-\delta t,\delta t)$ at time $t$ grows exponentially with a positive rate $\lambda_{1}-\sqrt{2\lambda_{1}}\delta$,
while for $\delta>\sqrt{\lambda_{1}/2}$, it converges to $0$.
 In the latter case,
 Proposition \ref{prop1} below shows that the upper bound of the
 mass
 outside $(-\delta t,\delta t)$ decreases exponentially with a negative rate.
The version of Theorem \ref{them1} has been
 proved recently for branching Brownian motions with branching rate given by a compactly supported measure in \cite{BH,Shiozawa1,Shiozawa2}.
 The idea of our proof is similar to that of \cite{Shiozawa2}: The upper bound and the lower bound of $\mathcal{X}^{\delta t}_{t}$ are considered separately and the proofs for convergence follow two main steps. The first step is to obtain the limit along lattice times. This is done via a Borel-Cantelli argument and thus requires the asymptotics of the expectation of $\mathcal{X}^{\delta t}_{t}$ (Lemma \ref{lem1} below). The second step is to extend the limits to all times. For the aforementioned class of branching Brownian motions, a key fact used in the proofs is that the particles alive at time $t$ located in $(-\delta t,\delta t)^{c}$ are the children of the particles alive at time $\lfloor t\rfloor$. However, this kind of property fails for the super-Brownian motions. We overcome this difficulty by appealing to a stochastic integral representation of super-Brownian motions (eq. \eqref{eq:martingale representation} below). This representation enables us to decompose the super-Brownian motion in terms of martingale measures and hence providing useful structural properties of super-Brownian motions.
Let us mention that the result, which corresponds to Theorem \ref{them1} for the aforementioned class of branching Brownian motions,
implies that the supremum of the support of the process,
denoted by $R_{t}$, grows linearly with rate $\sqrt{\lambda_{1}/2}$ as $t\to+\infty$ a.s. on the survival event. It is further shown in \cite{NS} that
 $$R_{t}=\sqrt{\frac{\lambda_{1}}{2}}t+Y_{t},$$
 where the conditional distribution of $Y_{t}$ on the survival event is convergent. However, this property no longer holds for the super-Brownian motions. As we show in Remark \ref{rm:supremum} below, for the $(B_{t},\psi)$-superprocess, the conditional distributions of $R_{t}-\sqrt{\lambda_{1}/2}\,t$ are not even tight.

The growth order of $\mathcal{X}^{\delta t}_{t}$ undergoes the phase transition at $\delta=\sqrt{\lambda_{1}/2}$. We further obtain the limiting distributions of the super-Brownian motion
at the critical phase in Theorem \ref{them2} below.
For $\nu\in\mathcal{M}_{loc}(\R)$
and $x\in\R$, we use $\nu+x$ to denote the measure induced by the shift operator $y\mapsto x+y$, that is, $\int_{\R}f(y)(\nu+x)(dy)=\int_{\R}f(y+x)\nu(dy)$ for all $f\in\mathcal{B}^{+}(\R)$.

\begin{theorem}\label{them2}
For every $\mu\in\mf$, $((X_{t}\pm \sqrt{\lambda_{1}/2}\,t)_{t\ge 0},\pp_{\mu})$
converges
in distribution to
$W^{h}_{\infty}(X)\eta_{\pm}(dx)$, where $\eta_{\pm}(dx)$
are (non-random) measures on $\R$ defined by $\eta_{\pm}(dx)=C_{\pm}\e^{\pm\sqrt{2\lambda_{1}}x}dx$ with $C_{\pm}$ being defined by \eqref{def-C}.
\end{theorem}

For branching Markov processes, results of the type of Theorem \ref{them2} have been established in recent years for various models. See, e.g., \cite{ABBS,ABK,BKLMZ} for spatially-homogeneous branching Brownian motions, \cite{Aidekon,HS,Madaule} for branching random walks, \cite{RSZ2} for branching L\'{e}vy processes, \cite{BM,HRS} for multitype branching Brownian motions, and \cite{Bocharov,Nishimori} for spatially-inhomogeneous branching Brownian motions.
On contrast, there is much less work for superprocesses.
Very recently, Ren et al. \cite{RSZ} shows that super-Brownian motion with a spatially-independent branching mechanism
translated by a centered term converges in distribution.
Later, Ren et al. \cite{RYZ} represents the limiting process as the limit of a sequence of Poisson random measures in which each atom is decorated by an independent copy of an auxiliary measure. As far as the authors know, there are no references on the vague convergence for superprocesses with spatially-dependent branching mechanisms. To prove Theorem \ref{them2}, we appeal to the skeleton techniques for superprocesses. Intuitively, under suitable assumptions, for a given superprocess $(X_{t})_{t\ge 0}$ there exists a related branching Markov process $(Z_{t})_{t\ge 0}$, called the skeleton, such that at each fixed time $t\ge 0$, the law of $Z_{t}$ may be coupled to the law of $X_{t}$ in such a way that given $X_{t}$, $Z_{t}$ has the law of a Poisson point process with random intensity determined by $X_{t}$. We exploit this fact and carry the long time behavior from the skeleton to the superprocess. Our idea is partly inspired by \cite{RYZ} where the skeleton techniques have been used successfully to establish the limiting distribution for super-Brownian motions with spatially-independent branching mechanisms.

Theorem \ref{them2} yields the following result on the convergence of the mass at the critical phase.

\begin{theorem}\label{them3}
For $\delta=\sqrt{\lambda_{1}/2}$ and $\mu\in\mc$, $(\mathcal{X}^{\delta t}_{t},\pp_{\mu})$ converges in distribution to $\frac{1}{\sqrt{2\lambda_{1}}}(C_{+}+C_{-})W^{h}_{\infty}(X)$.
\end{theorem}

The rest of this paper is organized as follows. In Section \ref{sec2} we derive the long time asymptotic properties of Feynman-Kac functionals related to the first and second moments of
superprocesses.
Section \ref{sec3} is devoted to the proof of Theorem \ref{them1}.
The proofs of Theorems \ref{them2} and \ref{them3} are given in Section \ref{sec4}.

\section{Estimates on the Feynman-Kac functionals}\label{sec2}

In this section, we show two lemmas related to the Feynman-Kac functionals of Brownian motions,
which will be used in the proofs of the main results.

Let $a(t)$ be a function on $[0,+\infty)$ with $a(t)=o(t)$ as $t\to+\infty$. For $\delta>0$,
define $R(t):=\delta t+a(t)$.
Let $A$ be a Borel set of $\R$ with $\inf A>-\infty$. Let
$b:[0,+\infty)\to [0,+\infty)$
be a function with $b(t)=o(t)$ as $t\to+\infty$.
For $r\in\R$ and $\Theta\subseteq\{\pm 1\}$, define $C_{\Theta}(r,A):=\{x\in\R:\ \theta x\in r+A\mbox{ for some }\theta\in\Theta\}$.

\begin{lemma}\label{lem:estimate}
Suppose $\delta\in (0,\sqrt{2\lambda_{1}})$.
\begin{description}
\item{(i)} For any $a\in \left(0,1-\frac{\delta}{\sqrt{2\lambda_{1}}}\right)$, there exist constants ${C_{1}},T_{1}>0$ such that for $t\ge  T_{1}$, $s\in [0,at]$ and $|x|\le b(t)$,
\begin{equation}\label{esti1}
\Big|\Bp_{x}\left[e_{\beta}(t-s),B_{t-s}\in C_{\Theta}(R(t),A)\right]-\e^{\lambda_{1}(t-s)}h(x)\int_{C_{\Theta}(R(t),A)}h(y)dy\Big|\le \e^{-{C_{1}}t}\e^{\lambda_{1}(t-s)-\sqrt{2\lambda_{1}}R(t)}.
\end{equation}
\item{(ii)} There exist constants ${C_{2}}>0$ and $T_{2}>1$ such that for $t\ge T_{2}$, $s\in [0,t-1]$ and $|x|\le b(t)$,
\begin{equation}\label{esti2}
\Bp_{x}\left[\int_{0}^{t-s}\gamma(B_{r})e_{\beta}(r)\Bp_{B_{r}}\big[e_{\beta}(t-s-r),B_{t-s-r}\in C_{\Theta}(R(t),A)\big]^{2}dr\right]\le {C_{2}}\e^{2\lambda_{1}(t-s)-2\sqrt{2\lambda_{1}}R(t)},
\end{equation}
where $\gamma$ is defined by \eqref{def-gamma}.
\end{description}
\end{lemma}

We remark here that for the special case where $\Theta=\{\pm 1\}$, $A=(0,+\infty)$ (correspondingly $C_{\Theta}(R(t),A)=\{y\in\R:\ |y|>R(t)\}$) and
$b(t)=b$ for some constant $b>0$,
the above two inequalities follow, respectively, from Lemma 3.8 and Lemma 3.9 of \cite{NS}.
Here we show the results for more general case where $A$ can be any left-bounded Borel set and
$b(t)=o(t)$.
Our proofs are based on
\cite[Section 3.3]{NS}.

\medskip

\noindent\textbf{Proof of Lemma \ref{lem:estimate}:}
(i) Let $p^{\beta}(t,x,y)$ and $p(t,x,y)$ be the transition densities of $P^{\beta}_{t}$ and $P_{t}$ respectively.
Let $q_{t}(x,y):=p^{\beta}(t,x,y)-p(t,x,y)-\e^{\lambda_{1}t}h(x)h(y)$.
We have
\begin{eqnarray}
&&\Bp_{x}\left[e_{\beta}(t-s),B_{t-s}\in C_{\Theta}(R(t),A)\right]
-\e^{\lambda_{1}(t-s)}h(x)\int_{C_{\Theta}(R(t),A)}h(y)dy\nonumber\\
&=&\Bp_{x}\left(B_{t-s}\in C_{\Theta}(R(t),A)\right)
+
\int_{C_{\Theta}(R(t),A)}q_{t-s}(x,y)dy.\nonumber
\end{eqnarray}
We note that for $R>0$ large enough such that $R+\inf A>0$,
$C_{\Theta}(R,A)\subseteq \{y\in\R:\ |y|\ge R+\inf A\}$. Thus for $t$ sufficiently large such that $R(t)-b(t)+\inf A>0$ and $|x|\le b(t)$,
\begin{equation}\label{esti3}
\Bp_{x}\left(B_{t-s}\in C_{\Theta}(R(t),A)\right)=\Bp_{0}\left(B_{t-s}+x\in C_{\Theta}(R(t),A)\right)\le \Bp_{0}\left(|B_{t-s}|\ge R(t)-b(t)-\inf A\right).
\end{equation}
On the other hand, it follows similarly as \cite[equation (3.19)]{NS} that for any $t\ge 1$ and $x\in\R$,
\begin{eqnarray}
\int_{C_{\Theta}(R,A)}q_{t}(x,y)dy&=&\int_{0}^{1}\Big(\int_{\R}p^{\beta}_{s}(x,z)\Bp_{z}\big(B_{t-s}\in C_{\Theta}(R,A)\Big)\beta(z)dz\big)ds\nonumber\\
&&+\int_{1}^{t}\Big[\int_{\R}\big(p^{\beta}_{s}(x,z)-\e^{\lambda_{1}s}h(x)h(z)\big)\Bp_{z}\big(B_{t-s}\in C_{\Theta}(R,A)\big)\beta(z)dz\Big]ds\nonumber\\
&&-\e^{\lambda_{1}t}h(x)\int_{t-1}^{+\infty}\e^{-\lambda_{1}s}\Big(\int_{\R}h(z)\beta(z)\Bp_{z}\left(B_{s}\in C_{\Theta}(R,A)\right)dz\Big)ds.\nonumber
\end{eqnarray}
Thus we have
$$\left|\int_{C_{\Theta}(R,A)}q_{t}(x,y)dy\right|\le (I)+(II)+(III),$$
where
\begin{eqnarray*}
(I)&=&\int_{0}^{1}\Big(\int_{\R}p^{\beta}_{s}(x,z)\Bp_{z}\big(|B_{t-s}|\ge R+\inf A\big)\beta(z)dz\Big)ds,\\
(II)&=&\int_{1}^{t}\Big[\int_{\R}\big(p^{\beta}_{s}(x,z)-\e^{\lambda_{1}s}h(x)h(z)\big)\Bp_{z}\big(|B_{t-s}|\ge R+\inf A\big)\beta(z)dz\Big]ds,\\
(III)&=&\e^{\lambda_{1}t}h(x)\int_{t-1}^{+\infty}\e^{-\lambda_{1}s}\Big(\int_{\R}h(z)\beta(z)\Bp_{z}\left(|B_{s}|\ge  R+\inf A\right)dz\Big)ds.
\end{eqnarray*}
The upper bounds for $(I),(II),(III)$ are established through Lemmas 3.5-3.7 of \cite{NS}.
These yield that
if supp$\beta\subset [-k,k]$ for some $k\in (0,+\infty)$, then there exist constants $c,C>0$ such that for all $x\in\R$, $t\ge 1$ and $R+\inf A>2k$,
\begin{eqnarray}
\left|\int_{C_{\Theta}(R,A)}q_{t}(x,y)dy\right|&\le &C\Big[h(x)\Bp_{0}\left(|B_{t}|>R+\inf A-k\right)\nonumber\\
&&+I_{c}(t,R+\inf A)+h(x)J(t,R+\inf A)\Big].\label{esti4}
\end{eqnarray}
Here $I_{c}$ and $J$ are defined by (3.15) and (3.16) of \cite{NS} respectively.
Using \eqref{esti3} and \eqref{esti4}, one can apply similar argument of \cite[Lemma 3.8]{NS} to prove \eqref{esti1}. We omit the details here.

\medskip

(ii) Noting that for $t$ large enough such that
$R(t)+\inf A\ge 0$, $C_{\Theta}(R(t),A)\subseteq \{y\in \R:\ |y|\ge R(t)+\inf A\}$, we have
\begin{eqnarray}
&&\Bp_{x}\left[\int_{0}^{t-s}\gamma(B_{r})e_{\beta}(r)\Bp_{B_{r}}\big[e_{\beta}(t-s-r),B_{t-s-r}\in C_{\Theta}(R(t),A)\big]^{2}dr\right]\nonumber\\
&\le&
\Bp_{x}\left[\int_{0}^{t-s}\gamma(B_{r})e_{\beta}(r)\Bp_{B_{r}}\big[e_{\beta}(t-s-r),|B_{t-s-r}|\ge R(t)+\inf A\big]^{2}dr\right].
\end{eqnarray}
Using the argument of \cite[Lemma 3.9]{NS} with minor modifications, one can prove \eqref{def-gamma}.
We omit the details.\qed

\bigskip

\begin{lemma}\label{lem:esti2}
Suppose the assumptions of Lemma \ref{lem:estimate}(i) hold. Then there exist $T>0$ and $\theta_{\pm}(t)$ such that for $t\ge T$, $s\in [0,at]$ and $|x|\le b(t)$,
\begin{equation}\label{esti2.1}
\theta_{-}(t)\le \frac{\Bp_{x}\big[e_{\beta}(t-s),B_{t-s}\in C_{\Theta}(R(t),A)\big]}{C_{\Theta}\left(\int_{A}\e^{-\sqrt{2\lambda_{1}}y}dy\right) h(x)\e^{\lambda_{1}(t-s)-\sqrt{2\lambda_{1}}R(t)}}\le\theta_{+}(t),
\end{equation}
where $\theta_{\pm}(t)\to 1$ as $t\to+\infty$ and $C_{\Theta}=C_{-},\ C_{+}$ and $(C_{+}+C_{-})$ accordingly as $\Theta=\{1\}$, $\{-1\}$ and $\{\pm 1\}$.
\end{lemma}

\proof Without loss of generality we assume in addition that $b(t)\to+\infty$ as $t\to+\infty$. Noting \eqref{eq:h}, we have
$$\int_{C_{\Theta}(R(t),A)}h(y)dy=\int_{C_{\Theta}(R(t),A)}\left(\int_{-\infty}^{+\infty}G_{\lambda_{1}}(y,z)\beta(z)h(z)dz\right)dy.$$
Using \eqref{estiforG} and the fact that $\beta$ is compactly supported, one can easily show by elementary calculation that
$$\int_{C_{\Theta}(R(t),A)}h(y)dy\sim C_{\Theta}\eta(A)\e^{-\sqrt{2\lambda_{1}}R(t)}\quad\mbox{ as }t\to+\infty,$$
where $\eta(A)=\int_{A}\e^{-\sqrt{2\lambda_{1}}y}dy$.
It then follows from Lemma \ref{lem:estimate}(i) that there exist constants ${c_{1}},T_{1}>0$ such that for $t\ge  T_{1}$, $s\in [0,at]$ and $|x|\le b(t)$,
\begin{equation}
\left|\frac{\Bp_{x}\big[e_{\beta}(t-s),B_{t-s}\in C_{\Theta}(R(t)+A)\big]}{C_{\Theta}\eta(A)h(x)\e^{\lambda_{1}(t-s)-\sqrt{2\lambda_{1}}R(t)}}-\frac{\int_{C_{\Theta}(R(t),A)}h(y)dy}{C_{\Theta}\eta(A)\e^{-\sqrt{2\lambda_{1}}R(t)}}\right|\le \frac{\e^{-{c_{1}}t}}{C_{\Theta}\eta(A)h(x)}.
\end{equation}
By \eqref{estimateforh}, there is a constant $c_{2}>0$ such that $h(x)\ge c_{2}\e^{-\sqrt{2\lambda_{1}}|x|}$ for all $x\in\R$.
So one has $\inf_{|x|\le b(t)}h(x)\ge c_{2}\e^{-\sqrt{2\lambda_{1}}b(t)}$.
Thus
$$\frac{\e^{-{c_{1}}t}}{C_{\Theta}\eta(A)h(x)}\le c_{3}\e^{-{c_{1}}t+\sqrt{2\lambda_{1}}b(t)}\to  0\quad\mbox{ as }t\to+\infty.$$
Hence we obtain  \eqref{esti2.1}
by setting $\theta_{\pm}(t)=\frac{\int_{C_{\Theta}(R(t),A)}h(y)dy}{C_{\Theta}\eta(A)\e^{-\sqrt{2\lambda_{1}}R(t)}}\pm c_{3}\e^{-{c_{1}}t+\sqrt{2\lambda_{1}}b(t)}$.\qed

\section{Proof of Theorem \ref{them1}}\label{sec3}

\subsection{Estimates on the first moment}

Put
$$\pi^{R}_{t}(x):=\Bp_{x}\left[e_{\beta}(t);|B_{t}|\ge R\right],\quad t\ge 0,x\in\R.$$
In this section we derive some estimates for $\pi^{R}_{t}(x)$, which will be used in the proof of Theorem \ref{them1}.

For any $\delta\ge 0$, we define
\begin{equation} \label{def:Lambda}
\Lambda_{\delta}:=\begin{cases}
         -\lambda_{1}+\sqrt{2\lambda_{1}}\delta
         \quad &\hbox{if } 0\le \delta < \sqrt{2\lambda_{1}},  \\
         \frac{\delta^{2}}{2}
         \quad &\hbox{if } \delta \ge \sqrt{2\lambda_{1}}.
        \end{cases}
\end{equation}
Obviously, $\Lambda_{\delta}<0$, $\Lambda_{\delta}=0$ and $\Lambda_{\delta}>0$ accordingly as $0\le \delta<\sqrt{\lambda_{1}/2}$, $\delta=\sqrt{\lambda_{1}/2}$ and $\delta>\sqrt{\lambda_{1}/2}$.

\begin{lemma}\label{lem1}
Suppose $\delta>0$. For any compact set $K\subset \R$,
\begin{equation}\label{lem1.eq2}
\lim_{t\to+\infty}\sup_{y\in K}\frac{\log \pi^{\delta t}_{t}(y)}{t}=\lim_{t\to+\infty}\inf_{y\in K}\frac{\log \pi^{\delta t}_{t}(x)}{t}=-\Lambda_{\delta}.
\end{equation}
\end{lemma}

\proof For $\delta\ge \sqrt{2\lambda_{1}}$, \eqref{lem1.eq2} is proved by \cite[Lemmas 4.4-4.5]{Shiozawa2}.
For $0<\delta<\sqrt{2\lambda_{1}}$, noting that $\pi^{\delta t}_{t}(x)=\Bp_{x}\left[e_{\beta}(t),B_{t}\in C_{\{\pm 1\}}(\delta t,(0,+\infty))\right]$, we get by Lemma \ref{lem:esti2} that for every compact set $K$, when $t$ is sufficiently large,
$$\theta_{-}(t)c_{1}\e^{\lambda_{1}t-\sqrt{2\lambda_{1}}\delta t}\le \pi^{\delta t}_{t}(y)\le \theta_{+}(t)c_{1}\e^{\lambda_{1}t-\sqrt{2\lambda_{1}}\delta t}\quad\forall y\in K,$$
where $c_{1}=(C_{+}+C_{-})\int_{0}^{+\infty}\e^{-\sqrt{2\lambda_{1}}y}dy$ and $\theta_{\pm}(t)\to 1$ as $t\to+\infty$. Thus \eqref{lem1.eq2} follows immediately.\qed

\bigskip

\begin{remark}\rm
We remark here that for any compact set $K\subset \R$ and $x\in\R$,
\begin{equation}\label{lem1.eq3}
\lim_{t\to+\infty}\sup_{y\in K}\frac{\log \pi^{0}_{t}(y)}{t}=\lim_{t\to+\infty}\frac{\log \pi^{0}_{t}(x)}{t}=\lambda_{1}.
\end{equation}
The second equality follows immediately by \cite[Theorem A.2]{Shiozawa2}. We shall show the first equality. Let $\epsilon>0$. It follows by \cite[Lemma 4.3]{Shiozawa2} that there is $p^{*}>1$ such that for all $p\in (1,p^{*})$,
$$c_{1}(p):=\sup_{y\in\R }\Bp_{x}\left[\sup_{t\ge 0}\mathrm{e}^{-p(\lambda_{1}+\epsilon)t}e_{p\beta}(t)\right]<+\infty.$$
By this and Jensen's inequality, we have
\begin{equation*}
\pi^{0}_{t}(y)=\Bp_{y}\left[e_{\beta}(t)\right]
\le\Bp_{y}\left[\e_{p\beta}(t)\right]^{1/p}=\mathrm{e}^{(\lambda_{1}+\epsilon)t}\Bp_{y}\left[\mathrm{e}^{-p(\lambda_{1}+\epsilon)t}e_{p\beta}(t)\right]
\le c_{1}(p)^{1/p}\e^{(\lambda_{1}+\epsilon)t}
\end{equation*}
for every $t\ge 0$ and $y\in\R $. Thus
$$\limsup_{t\to+\infty}\sup_{y\in K }\frac{\log \pi^{0}_{t}(y)}{t}\le \lambda_{1}+\epsilon.$$
This implies the first identity of \eqref{lem1.eq3}.
\end{remark}

\bigskip

\begin{lemma}\label{lem2}
\begin{description}
\item{(i)} For every $\sigma>0$, there exists a
constant
$C_{3}=C_{3}(\sigma)>0$,
such that for any $0<\theta<\delta<+\infty$,
$$\pi^{\delta t}_{s}(x)\le
C_{3}\cdot
 (\theta t)^{-1}\e^{-\frac{\theta^{2}t^{2}}{2\sigma}},\quad\forall s\in (0,\sigma],\ |x|\le(\delta-\theta)t,$$
when $t$ is sufficiently large.

\item{(ii)} For every $\delta\ge 0$ and $\sigma>0$, there exists a constant
$C_{4}=C_{4}(\delta,\sigma)>0$
such that for any $s\in (0,\sigma]$, $t\ge s$ and $|x|\ge \delta(t-s)$,
$$\pi^{\delta t}_{s}(x)\ge {C_{4}}.$$

\item{(iii)}
If $\delta>\sqrt{\lambda_{1}/2}$, then
$x\mapsto \int_{0}^{+\infty}\pi^{\delta s}_{s}(x)ds$ is a locally bounded function on $\R $.
\end{description}
\end{lemma}

\proof (i) Note that for every $s\ge 0$,
$$e_{\beta}(s)=1+\int_{0}^{s}e_{\beta}(r)\beta(B_{r})dr.$$
We have
\begin{eqnarray}
\pi^{\delta t}_{s}(x)&=&\Bp_{x}\left[e_{\beta}(s);|B_{s}|\ge \delta t\right]\nonumber\\
&=&\Bp_{x}\left(|B_{s}|\ge \delta t\right)+\Bp_{x}\left[\int_{0}^{s}e_{\beta}(r)\beta(B_{r})\1_{\{|B_{s}|\ge \delta t\}}dr\right]\nonumber\\
&=:&\mathrm{I}(x,s,t)+\mathrm{II}(x,s,t).\notag
\end{eqnarray}
For $R\ge 0$, let
$$G(R):=\Bp_{0}\left(|B_{1}|\ge R\right)=\sqrt{\frac{2}{\pi}}\int_{R}^{+\infty}\e^{-\frac{y^{2}}{2}}dy.$$
Then for $s\in (0,\sigma]$ and $|x|\le (\delta-\theta)t$,
\begin{equation}
\mathrm{I}(x,s,t)=\Bp_{0}\left(|B_{s}+x|\ge \delta t\right)\le\Bp_{0}\left(|B_{s}|\ge \delta t-|x|\right)\le\Bp_{0}\left(|B_{1}|\ge \frac{\theta t}{\sqrt{\sigma}}\right)=G\left(\frac{\theta t}{\sqrt{\sigma}}\right).\label{lem2.1}
\end{equation}
Suppose $\mbox{supp}\beta\subset[-k,k]$ for some $k\in(0,+\infty)$.
By Markov property, we have for $x\in\R$ and $s\in (0,\sigma]$,
\begin{align}
\mathrm{II}(x,s,t)&=\Bp_{x}\left[\int_{0}^{s}e_{\beta}(r)\beta(B_{r})\Bp_{B_{r}}\left(|B_{s-r}|\ge \delta t\right)dr\right]\notag\\
&=\Bp_{x}\left[\int_{0}^{s}e_{\beta}(r)\beta(B_{r})1_{\{|B_{r}|\le k\}}\Bp_{B_{r}}\left(|B_{s-r}|\ge \delta t\right)dr\right]\notag\\
&\le\Bp_{x}\left[\int_{0}^{s}e_{\beta}(r)\beta^{+}(B_{r})\Bp_{0}\left(|B_{s-r}|\ge \delta t-k\right)dr\right]\notag\\
&\le \int_{0}^{\sigma}\e^{\|\beta^{+}\|_{\infty}r}\|\beta^{+}\|_{\infty}G\left(\frac{\delta t-k}{\sqrt{\sigma-r}}\right)dr\notag\\
&\le c_{1}G\left(\frac{\delta t-k}{\sqrt{\sigma}}\right)\label{lem2.2}
\end{align}
for some $c_{1}=c_{1}(\sigma)>0$. Note that for $t\ge k/(\delta-\theta)$,
$G((\delta t-k)/\sqrt{\sigma})\le G(\theta t/\sqrt{\sigma})$.
It follows from
 \eqref{lem2.1} and \eqref{lem2.2} that for $t\ge k/(\delta-\theta)$, $|x|\le (\delta-\theta)t$ and $s\in (0,\sigma]$,
$$\pi^{\delta t}_{s}(x)\le (1+c_{1})G\left(\frac{\theta t}{\sqrt{\sigma}}\right).$$
Thus (i) follows by the fact that $G(R)\sim \sqrt{2/\pi}R^{-1}\e^{-R^{2}/2}$ as $R\to+\infty$.

\medskip

(ii) We have
$$\pi^{\delta t}_{s}(x)=\Bp_{x}\left[e_{\beta}(s);|B_{s}|\ge \delta t\right]\ge \mathrm{e}^{-\|\beta^{-}\|_{\infty}s}\Bp_{x}\left(|B_{s}|\ge\delta t\right).$$
Note that for $x\ge \delta(t-s)$,
$$\Bp_{x}\left(|B_{s}|\ge \delta t\right)\ge \Bp_{x}\left(B_{s}\ge \delta t\right)=\Bp_{0}\left(B_{s}\ge \delta t-x\right)\ge \Bp_{0}\left(B_{s}\ge \delta s\right)=\Bp_{0}\left(B_{1}\ge\delta\sqrt{s}\right).$$
Similarly one can show that $\Bp_{x}\left(|B_{s}|\ge \delta t\right)\ge \Bp_{0}\left(B_{1}\le -\delta\sqrt{s}\right)$ for $x\le -\delta(t-s)$.
Hence we get (ii) by setting $C_{4}=\mathrm{e}^{-\|\beta^{-}\|_{\infty}\sigma}\Bp_{0}\left(B_{1}\ge \delta\sqrt{\sigma}\right)$.

\medskip

(iii) By Lemma \ref{lem1}, for $\delta>\sqrt{\lambda_{1}/2}$ and any compact set $K\subset \R$,
$$\lim_{t\to+\infty}\sup_{x\in K}\frac{\log \pi^{\delta t}_{t}(x)}{t}=-\Lambda_{\delta}<0.$$
So there is some $T>0$ such that for all $s\ge T$,
$\sup_{x\in K}\pi^{\delta s}_{s}(x)\le \exp\{-\Lambda_{\delta}s/{2}\}$.
We also note that $\pi^{\delta s}_{s}(x)=\Bp_{x}\left[e_{\beta}(s);
|B_{s}|\ge \delta s\right]\le \e^{\|\beta^{+}\|_{\infty}s}$ for all $x\in K$ and $s\ge 0$. Hence
$$\sup_{x\in K}\int_{0}^{+\infty}\pi^{\delta s}_{s}(x)ds\le \int_{0}^{T}\e^{\|\beta^{+}\|_{\infty}s}ds+\int_{T}^{+\infty}\e^{-\Lambda_{\delta}s/2}ds<+\infty.$$\qed

\subsection{The upper bound of $\mathcal{X}^{\delta t}_{t}$}

For $t\ge 0$, let $\mathcal{F}_{t}$ denote the $\sigma$-field generated by $\{X_{s}:0\le s\le t\}$.
It follows immediately from Lemma \ref{lem2}(ii) that for any $n\in\mathbb{N}$ and $t\in [n\sigma, (n+1)\sigma)$,
$$\pi^{\delta (n+1)\sigma}_{(n+1)\sigma-t}(x)\ge {C_{4}}, \quad\forall |x|\ge \delta t.$$
Thus for $t\in [n\sigma,(n+1)\sigma)$ and $\mu\in\mf$,
\begin{equation}\label{fact:6}
\mathcal{X}^{\delta t}_{t}=\langle \1_{(-\delta t,\delta t)^{c}},X_{t}\rangle \le C^{-1}_{4}\langle \pi^{\delta (n+1)\sigma}_{(n+1)\sigma-t},X_{t}\rangle=C^{-1}_{4}\pp_{\mu}\left[\mathcal{X}^{\delta(n+1)\sigma}_{(n+1)\sigma}|\mathcal{F}_{t}\right]\quad\pp_{\mu}\mbox{-a.s.}
\end{equation}
Here the last equality follows by the Markov property of $X_{t}$.
Hence to get an upper bound for $\mathcal{X}^{\delta t}_{t}$ we only need to
compute
$\pp_{\mu}\left[\mathcal{X}^{\delta (n+1)\sigma}_{(n+1)\sigma}\,|\,\mathcal{F}_{t}\right]$.

\begin{lemma}\label{lem0}
For any $\delta\ge 0$, $\sigma>0$ and $\mu\in\mc$,
$$\limsup_{n\to+\infty}\frac{\log \mathcal{X}^{\delta n\sigma}_{n\sigma}}{n\sigma}\le -\Lambda_{\delta}\quad \pp_{\mu}\text{-a.s.}$$
\end{lemma}

\proof Let $\epsilon>0$. We have
\begin{align*}
\pp_{\mu}\left(\frac{\log \mathcal{X}^{\delta n\sigma}_{n\sigma}}{n\sigma}\ge -\Lambda_{\delta}+\epsilon\right)
&=\pp_{\mu}\left(\mathcal{X}^{\delta n \sigma}_{n\sigma}\ge \e^{(-\Lambda_{\delta}+\epsilon)n\sigma}\right)\\
&\le \e^{(\Lambda_{\delta}-\epsilon)n\sigma}\langle \pi^{\delta n \sigma}_{n\sigma},\mu\rangle\\
&=\e^{n\sigma\left(\frac{\log \langle \pi^{\delta n \sigma}_{n\sigma},\mu\rangle}{n\sigma}+\Lambda_{\delta}-\epsilon\right)}.
\end{align*}
Since $\mu$ is compactly supported, Lemma \ref{lem1} and \eqref{lem1.eq3} imply that
$$\limsup_{n\to +\infty}\frac{\log \langle \pi^{\delta n \sigma}_{n\sigma},\mu\rangle}{n\sigma}\le -\Lambda_{\delta}.$$
Thus when $n$ is large enough, we have
$$\pp_{\mu}\left(\frac{\log \mathcal{X}^{\delta n\sigma}_{n\sigma}}{n\sigma}\ge -\Lambda_{\delta}+\epsilon\right)\le \e^{-\epsilon n\sigma},$$
which in turn implies that
$\sum_{n=0}^{+\infty}\pp_{\mu}\left(\frac{\log \mathcal{X}^{\delta n\sigma}_{n\sigma}}{n\sigma}\ge -\Lambda_{\delta}+\epsilon\right)<+\infty$. Hence this lemma follows immediately by Borel-Cantelli lemma.\qed

\bigskip

It follows from \cite[Corollary 2.18]{F1992} that for any $f\in\B_{b}(\R )$, $t\ge 0$ and $\mu\in\me$,
\begin{equation}\label{eq:martingale representation}
\langle f,X_{t}\rangle =\langle P^{\beta}_{t}f,X_{0}\rangle +\int_{0}^{t}\int_{-\infty}^{+\infty }P^{\beta}_{t-s}f(x)M(ds,dx)\quad \pp_{\mu}\mbox{-a.s.}
\end{equation}
where for every $T>0$, $[0,T]\ni t\mapsto \int_{0}^{t}\int_{-\infty}^{+\infty }P^{\beta}_{t-s}f(x)M(ds,dx)$ is a square-integrable $\mathcal{F}_{t}$-martingale with quadratic variation
$t\mapsto \int_{0}^{t}\langle 2\gamma (P^{\beta}_{t-s}f)^{2},X_{s}\rangle ds$.

\begin{lemma}\label{lem4}
Suppose $\delta\ge 0$, $\sigma>0$ and $\mu\in\mc$. Then for any $\epsilon>0$,
$$\lim_{n\to+\infty}\e^{\left(\frac{1}{2}\Lambda_{\delta}-\epsilon\right)(n+1)\sigma}\sup_{t\in [n\sigma,(n+1)\sigma]}\left|
\pp_{\mu}\left[\mathcal{X}^{\delta(n+1)\sigma}_{(n+1)\sigma}|\mathcal{F}_{t}\right]-\pp_{\mu}\left[\mathcal{X}^{\delta(n+1)\sigma}_{(n+1)\sigma}|\mathcal{F}_{n\sigma}\right]\right|=0\quad \pp_{\mu}\text{-a.s.}$$
\end{lemma}

\proof By \eqref{eq:martingale representation}, we have
\begin{equation}\nonumber
\mathcal{X}^{\delta(n+1)\sigma}_{(n+1)\sigma}=\langle \pi^{\delta(n+1)\sigma}_{(n+1)\sigma},X_{0}\rangle +\int_{0}^{(n+1)\sigma}\int_{-\infty}^{+\infty }\pi^{\delta(n+1)\sigma}_{(n+1)\sigma-s}(x)M(ds,dx),
\end{equation}
where $[0,(n+1)\sigma]\ni t\mapsto \int_{0}^{t}\int_{-\infty}^{+\infty }\pi^{\delta(n+1)\sigma}_{(n+1)\sigma-s}(x)M(ds,dx)$
is a square-integrable $\pp_{\mu}$-martingale with quadratic variation
$t\mapsto \int_{0}^{t}\langle 2\gamma (\pi^{\delta(n+1)\sigma}_{(n+1)\sigma-s})^{2},X_{s}\rangle ds$. Thus
\begin{align*}
&\pp_{\mu}\left[\mathcal{X}^{\delta(n+1)\sigma}_{(n+1)\sigma}|\mathcal{F}_{t}\right]-\pp_{\mu}\left[\mathcal{X}^{\delta(n+1)\sigma}_{(n+1)\sigma}|\mathcal{F}_{n\sigma}\right]\\
&=\pp_{\mu}\left[\int_{0}^{(n+1)\sigma}\int_{-\infty}^{+\infty }\pi^{\delta(n+1)\sigma}_{(n+1)\sigma-s}(x)M(ds,dx)\Big|\mathcal{F}_{t}\right]
-\pp_{\mu}\left[\int_{0}^{(n+1)\sigma}\int_{-\infty}^{+\infty }\pi^{\delta(n+1)\sigma}_{(n+1)\sigma-s}(x)M(ds,dx)\Big|\mathcal{F}_{n\sigma}\right]\\
&=\int_{n\sigma}^{t}\int_{-\infty}^{+\infty }\pi^{\delta(n+1)\sigma}_{(n+1)\sigma-s}(x)M(ds,dx).
\end{align*}
By this and the $L^{2}$-maximum inequality for martingales, we have
\begin{align}
&\pp_{\mu}\left[\sup_{t\in [n\sigma,(n+1)\sigma]}\left|\pp_{\mu}\left[\mathcal{X}^{\delta(n+1)\sigma}_{(n+1)\sigma}|\mathcal{F}_{t}\right]-\pp_{\mu}\left[\mathcal{X}^{\delta(n+1)\sigma}_{(n+1)\sigma}|\mathcal{F}_{n\sigma}\right]\right|^{2}\right]\notag\\
&\le 4\pp_{\mu}\left[\left(\int_{n\sigma}^{(n+1)\sigma}\int_{-\infty}^{+\infty }\pi^{\delta(n+1)\sigma}_{(n+1)\sigma-s}(x)M(ds,dx)\right)^{2}\right]\notag\\
&=4\pp_{\mu}\left[\int_{n\sigma}^{(n+1)\sigma}\langle 2\gamma\left(\pi^{\delta(n+1)\sigma}_{(n+1)\sigma-s}\right)^{2},X_{s}\rangle ds\right]\notag\\
&\le 8\|\gamma\|_{\infty}\e^{\|\beta^{+}\|_{\infty}\sigma}\pp_{\mu}\left[\int_{n\sigma}^{(n+1)\sigma}\langle \pi^{\delta(n+1)\sigma}_{(n+1)\sigma-s},X_{s}\rangle ds\right]\notag\\
&=8\|\gamma\|_{\infty}\e^{\|\beta^{+}\|_{\infty}\sigma}\sigma\langle \pi^{\delta(n+1)\sigma}_{(n+1)\sigma},\mu\rangle.\label{lem3.5.1}
\end{align}
The second inequality
is because $\pi^{\delta(n+1)\sigma}_{(n+1)\sigma-s}(x)\le \Bp_{x}\left[e_{\beta}((n+1)\sigma-s)\right]\le\mathrm{e}^{\|\beta^{+}\|_{\infty}((n+1)\sigma-s)}$.
Let $K$ be the compact support of $\mu$. Lemma \ref{lem1} implies that for any $\epsilon>0$, there is $T>0$ such that
$$\sup_{x\in K}\pi^{\delta t}_{t}(x)\le \e^{(-\Lambda_{\delta}+\epsilon)t}\quad \forall t\ge T.$$
Thus one has
$$\langle \pi^{\delta(n+1)\sigma}_{(n+1)\sigma},\mu\rangle\le \mathrm{e}^{(-\Lambda_{\delta}+\epsilon)(n+1)\sigma}\langle 1,\mu\rangle$$
for $n$ sufficiently large.
Putting this back to \eqref{lem3.5.1}, one gets
\begin{eqnarray*}
&&\pp_{\mu}\left[\left(\e^{(\frac{1}{2}\Lambda_{\delta}-\epsilon)(n+1)\sigma}\sup_{t\in [n\sigma,(n+1)\sigma]}\left|\pp_{\mu}\left[\mathcal{X}^{\delta(n+1)\sigma}_{(n+1)\sigma}|\mathcal{F}_{t}\right]-\pp_{\mu}\left[\mathcal{X}^{\delta(n+1)\sigma}_{(n+1)\sigma}|\mathcal{F}_{n\sigma}\right]
\right|\right)^{2}
\right]\\
&\le&c_{1}\e^{(\Lambda_{\delta}-2\epsilon)(n+1)\sigma}\langle\pi^{\delta(n+1)\sigma}_{(n+1)\sigma},\mu\rangle \le c_{1}\e^{-\epsilon (n+1)\sigma}\langle 1,\mu\rangle,
\end{eqnarray*}
for some constant $c_{1}=c_{1}(\sigma)>0$. This implies that
$$\sum_{n=0}^{+\infty}\pp_{\mu}\left[\left(\e^{(\frac{1}{2}\Lambda_{\delta}-\epsilon)(n+1)\sigma}\sup_{t\in [n\sigma,(n+1)\sigma]}\left|\pp_{\mu}\left[\mathcal{X}^{\delta(n+1)\sigma}_{(n+1)\sigma}|\mathcal{F}_{t}\right]-\pp_{\mu}\left[\mathcal{X}^{\delta(n+1)\sigma}_{(n+1)\sigma}|\mathcal{F}_{n\sigma}\right]
\right|\right)^{2}
\right]<+\infty.$$
The lemma follows by Borel-Cantelli lemma.\qed

\bigskip

\begin{lemma}\label{lem5}
Suppose $\sigma>0$ and $\mu\in\mc$.
There is a constant $C_{5}=C_{5}(\sigma)>0$
such that
for any
$0<\theta<\delta<+\infty$ and $\epsilon>0$ the following inequality holds $\pp_{\mu}$-a.s. for $n$ sufficiently large.
\begin{align*}
\pp_{\mu}\left[\mathcal{X}^{\delta(n+1)\sigma}_{(n+1)\sigma}|\mathcal{F}_{n\sigma}\right]
\le
C_{5}
\left[(\theta n \sigma)^{-1}\e^{-\frac{\theta^{2}}{2}n^{2}\sigma +(\lambda_{1}+\sqrt{2\lambda_{1}}
(\delta-\theta))n\sigma}W^{h}_{n\sigma}(X) +\e^{(-\Lambda_{\delta-\theta}+\epsilon)n\sigma}\right].
\end{align*}
\end{lemma}

\proof Fix arbitrary  $\sigma>0$ and $\mu\in\mc$. By Markov property, we have for any $0<\theta<\delta<+\infty$,
\begin{align*}
\pp_{\mu}\left[\mathcal{X}^{\delta(n+1)\sigma}_{(n+1)\sigma}|\mathcal{F}_{n\sigma}\right]
&=\langle \pi^{\delta(n+1)\sigma}_{\sigma},X_{n\sigma}\rangle\\
&=\langle \pi^{\delta(n+1)\sigma}_{\sigma}\1_{(-(\delta-\theta)n\sigma,(\delta-\theta)n\sigma)},X_{n\sigma}\rangle+\langle \pi^{\delta(n+1)\sigma}_{\sigma}\1_{(-(\delta-\theta)n\sigma,(\delta-\theta)n\sigma)^{c}},X_{n\sigma}\rangle\\
&=:\mathrm{I}(n,\delta,\theta)+\mathrm{II}(n,\delta,\theta).
\end{align*}
It follows from Lemma \ref{lem2}(i) that when $n$ is sufficiently large,
$$\pi^{\delta(n+1)\sigma}_{\sigma}(x)\le {C_{3}}(\theta(n+1)\sigma)^{-1}\e^{-\frac{\theta^{2}}{2}(n+1)^{2}\sigma},\quad \forall |x|<(\delta-\theta)n\sigma,$$
where $C_{3}=C_{3}(\sigma)>0$. Thus we have
\begin{align}
\mathrm{I}(n,\delta,\theta)&\le {C_{3}}(\theta(n+1)\sigma)^{-1}\e^{-\frac{\theta^{2}}{2}n^{2}\sigma }\langle \1_{(-(\delta-\theta)n\sigma,(\delta-\theta)n\sigma)},X_{n\sigma}\rangle\notag\\
&\le {C_{3}}(\theta(n+1)\sigma)^{-1}\e^{-\frac{\theta^{2}}{2}n^{2}\sigma }\e^{\lambda_{1}n\sigma}\Big\langle \e^{-\lambda_{1}n\sigma}
\frac{h}{\inf_{|x|<(\delta-\theta)n\sigma}h(x)},X_{n\sigma}\Big\rangle\notag\\
&={C_{3}}(\theta(n+1)\sigma)^{-1}\e^{-\frac{\theta^{2}}{2}n^{2}\sigma+\lambda_{1}n\sigma}\left(\inf_{|x|<(\delta-\theta)n\sigma}h(x)\right)^{-1}W^{h}_{n\sigma}(X).
\label{lem5.3}
\end{align}
The continuity of $h$ and \eqref{estimateforh} imply that $\inf_{|x|<(\delta-\theta)n\sigma}h(x)\ge c_{1}\e^{-\sqrt{2\lambda_{1}}(\delta-\theta)n\sigma}$ for $n$ sufficiently large. Thus we get by \eqref{lem5.3} that $\pp_{\mu}$-a.s.
\begin{equation}\label{lem5.1}
\mathrm{I}(n,\delta,\theta)\le c_{2}\big(\theta n\sigma\big)^{-1}\e^{-\frac{\theta^{2}}{2}n^{2}\sigma +(\lambda_{1}+\sqrt{2\lambda_{1}}(\delta-\theta))n\sigma}
W^{h}_{n\sigma}(X)
\end{equation}
for $n$ sufficiently large, where $c_{2}=c_{2}(\sigma)>0$.
On the other hand, by Lemma \ref{lem0} we have
$\limsup_{n\to+\infty}\log \mathcal{X}^{(\delta-\theta)n\sigma}_{n\sigma}/n\sigma\le -\Lambda_{\delta-\theta}$ $\pp_{\mu}$-a.s. Thus for any $\epsilon>0$,
$$\pp_{\mu}\left(\mathcal{X}^{(\delta-\theta)n\sigma}_{n\sigma}\le \e^{(-\Lambda_{\delta-\theta}+\epsilon)n\sigma}\text{ for $n$ sufficiently large}\right)=1.$$
Note that by definition $\mathrm{II}(n,\delta,\theta)\le \|\pi^{\delta(n+1)\sigma}_{\sigma}\|_{\infty}\mathcal{X}^{(\delta-\theta)n\sigma}_{n\sigma}\le \e^{\|\beta^{+}\|_{\infty}\sigma}\mathcal{X}^{(\delta-\theta)n\sigma}_{n\sigma}$ for every $n\in\mathbb{N}$. We get that
\begin{equation}\label{lem5.2}
\pp_{\mu}\left(\mathrm{II}(n,\delta,\theta)\le c_{3}e^{(-\Lambda_{\delta-\theta}+\epsilon)n\sigma}\text{ for $n$ sufficiently large}\right)=1
\end{equation}
for $c_{3}=\e^{\|\beta^{+}\|_{\infty}\sigma}$.
This lemma follows immediately by combining \eqref{lem5.1} and \eqref{lem5.2}.\qed

\bigskip

\begin{lemma}\label{lem6}
For any $\delta\ge 0$, $\sigma>0$ and $\mu\in\mc$,
$$\limsup_{n\to+\infty}\frac{\log \pp_{\mu}\left[\mathcal{X}^{\delta(n+1)\sigma}_{(n+1)\sigma}|\mathcal{F}_{n\sigma}\right]}{n\sigma}\le -\Lambda_{\delta}\quad \pp_{\mu}\mbox{-a.s.}$$
\end{lemma}

\proof First we consider $\delta>0$. It follows by Lemma \ref{lem5} that for any $0<\theta<\delta$ and $\epsilon>0$, $\pp_{\mu}$-a.s.
\begin{align}\label{lem6.1}
\e^{\Lambda_{\delta}n\sigma}\pp_{\mu}\left[\mathcal{X}^{\delta(n+1)\sigma}_{(n+1)\sigma}|\mathcal{F}_{n\sigma}\right]&
\le {C_{5}}(\sigma)\Big[\big(\theta n \sigma\big)^{-1}\e^{-\frac{\theta^{2}}{2}n^{2}\sigma+(\Lambda_{\delta}+\lambda_{1}+\sqrt{2\lambda_{1}}(\delta-\theta))n\sigma}W^{h}_{n\sigma}(X)\notag\\
&\quad +\e^{(\Lambda_{\delta}-\Lambda_{\delta-\theta}+\epsilon)n\sigma}\Big]
\end{align}
for $n$ sufficiently large.
Since $\delta\mapsto \Lambda_{\delta}$ is nondecreasing and continuous on $(0,+\infty)$ and that $\pp_{\mu}\left(W^{h}_{\infty}(X)<+\infty\right)=1$, one can choose $\theta$ so small that $\Lambda_{\delta}-\Lambda_{\delta-\theta}<\epsilon$.
We also note that for fixed $\delta$ and $\theta$,  the first term on the right hand side of \eqref{lem6.1}
 converges to $0$ $\pp_{\mu}$-a.s. as $n\to+\infty$. Thus we get by \eqref{lem6.1} that
$$\pp_{\mu}\left(\e^{\Lambda_{\delta}n\sigma}\pp_{\mu}\left[\mathcal{X}^{\delta(n+1)\sigma}_{(n+1)\sigma}|\mathcal{F}_{n\sigma}\right]\le 2\e^{2\epsilon n\sigma}
\text{ for $n$ sufficiently large}\right)=1.$$
Hence we prove this lemma for $\delta>0$. Now we suppose $\delta=0$. By Markov property, we have
$$\pp_{\mu}\left[\mathcal{X}^{0}_{(n+1)\sigma}|\mathcal{F}_{n\sigma}\right]
=\langle \pi^{0}_{\sigma},X_{n\sigma}\rangle\le \e^{\|\beta^{+}\|_{\infty}\sigma}\mathcal{X}^{0}_{n\sigma}.$$
It follows by Lemma \ref{lem0} that
$$\limsup_{n\to+\infty}\frac{\log \pp_{\mu}\left[\mathcal{X}^{0}_{(n+1)\sigma}|\mathcal{F}_{n\sigma}\right]}{n\sigma}\le
\limsup_{n\to+\infty}\frac{\log \mathcal{X}^{0}_{n\sigma}}{n\sigma}\le -\Lambda_{0}\quad \pp_{\mu}\mbox{-a.s.}$$
Hence we complete the proof.\qed

\bigskip

\begin{proposition}\label{prop1}
Suppose $\mu\in\mc$. For any $\delta>\sqrt{\lambda_{1}/2}$,
\begin{equation}\label{prop1.1}
\limsup_{t\to+\infty}\frac{\log \mathcal{X}^{\delta t}_{t}}{t}\le -\frac{1}{2}\Lambda_{\delta}\quad \pp_{\mu}\mbox{-a.s.,}
\end{equation}
and for any $0\le \delta\le \sqrt{\lambda_{1}/2}$,
\begin{equation}\label{prop1.2}
\limsup_{t\to +\infty}\frac{\log \mathcal{X}^{\delta t}_{t}}{t}\le -\Lambda_{\delta}\quad \pp_{\mu}\mbox{-a.s.}
\end{equation}
\end{proposition}

\proof Let $\sigma>0$. By \eqref{fact:6}, we have for any $n\in\mathbb{N}$ and $t\in [n\sigma,(n+1)\sigma)$,
\begin{equation}\label{prop1.3}
\mathcal{X}^{\delta_{t}}_{t}\le C^{-1}_{4}\pp_{\mu}\left[\mathcal{X}^{\delta(n+1)\sigma}_{(n+1)\sigma}|\mathcal{F}_{t}\right]\quad\pp_{\mu}\mbox{-a.s.}
\end{equation}
One can decompose $\pp_{\mu}\left[\mathcal{X}^{\delta(n+1)\sigma}_{(n+1)\sigma}|\mathcal{F}_{t}\right]$ as
$\mathrm{I}(n,\sigma,t)+\mathrm{II}(n,\sigma)$,
 where
 $$\mathrm{I}(n,\sigma,t):=\pp_{\mu}\left[\mathcal{X}^{\delta(n+1)\sigma}_{(n+1)\sigma}|\mathcal{F}_{t}\right]-\pp_{\mu}\left[\mathcal{X}^{\delta(n+1)\sigma}_{(n+1)\sigma}|
 \mathcal{F}_{n\sigma}\right]\mbox{ and }
 \mathrm{II}(n,\sigma):=\pp_{\mu}\left[\mathcal{X}^{\delta(n+1)\sigma}_{(n+1)\sigma}|\mathcal{F}_{n\sigma}\right].$$
 It follows by Lemma \ref{lem6} that for any $\epsilon>0$,
 \begin{equation}
 \pp_{\mu}\left(\mathrm{II}(n,\sigma)\le \e^{(-\Lambda_{\delta}+\epsilon)n\sigma}\text{ for $n$ sufficiently large}\right)=1.
 \end{equation}
 On the other hand, by Lemma \ref{lem4} we have
 \begin{equation}\label{prop1.4}
 \pp_{\mu}\left(\sup_{t\in [n\sigma,(n+1)\sigma)}\mathrm{I}(n,\sigma,t)\le \epsilon e^{(-\frac{1}{2}\Lambda_{\delta}+\epsilon)(n+1)\sigma}
 \text{ for $n$ sufficiently large}\right)=1
 \end{equation}
 Combining \eqref{prop1.3}-\eqref{prop1.4}, we get
 $$\pp_{\mu}\left(\sup_{t\in [n\sigma,(n+1)\sigma)}\mathcal{X}^{\delta t}_{t}\le C^{-1}_{4}\left(\e^{(-\Lambda_{\delta}+\epsilon)n\sigma}+\epsilon e^{(-\frac{1}{2}\Lambda_{\delta}+\epsilon)(n+1)\sigma}\right)\text{ for $n$ sufficiently large}\right)=1.$$
 It follows immediately that
 \begin{equation}\label{prop1.5}
 \limsup_{t\to+\infty}\frac{\log \mathcal{X}^{\delta t}_{t}}{t}\le (-\Lambda_{\delta})\vee \left(-\frac{1}{2}\Lambda_{\delta}\right)\quad \pp_{\mu}\mbox{-a.s.}
 \end{equation}
 If $0\le \delta\le \sqrt{\lambda_{1}/2}$, then $-\Lambda_{\delta}\ge -\Lambda_{\delta}/2\ge 0$, and  \eqref{prop1.2} follows directly from \eqref{prop1.5}.
Otherwise if $\delta>\sqrt{\lambda_{1}/2}$, then $-\Lambda_{\delta}<-\frac{1}{2}\Lambda_{\delta}<0$,
and hence \eqref{prop1.1} follows.\qed

\bigskip

\subsection{The lower bound of $\mathcal{X}^{\delta t}_{t}$}

Let $p^{\beta}(t,x,y)$ be the transition density of $P^{\beta}_{t}$. It is easy to see that
\begin{equation}\label{1.24}
\e^{-\|\beta^{-}\|_{\infty}t}p(t,x,y)\le p^{\beta}(t,x,y)\le \e^{\|\beta^{+}\|_{\infty}t}p(t,x,y)\quad\forall t\ge 0,\ x,y\in\R .
\end{equation}
Here $p(t,x,y)$ is the transition density of a Brownian motion on $\R $.
Let $P^{h}_{t}$ be the semigroup obtained from $P^{\beta}_{t}$
through Doob's $h$-transform, that is,
\begin{equation}\label{eq:pht}
P^{h}_{t}f(x)=\frac{\e^{-\lambda_{1}t}}{h(x)}P^{\beta}_{t}(hf)(x)\quad \forall t\ge 0,\ x\in\R ,\ f\in\B^{+}(\R ).
\end{equation}
Then $P^{h}_{t}$ has a transition density with respect to the measure $h^{2}(y)dy$, which is given by
\begin{equation}\label{eq:density pht}
p^{h}(t,x,y)=\frac{\e^{-\lambda_{1}t}p^{\beta}(t,x,y)}{h(x)h(y)}\quad \forall t\ge 0,\ x,y\in\R .
\end{equation}
It is proved in \cite{CRY} that there exists a constant $a>0$ such that
$$\left|p^{h}(t,x,y)-1\right|\le e^{-a(t-1)}p^{h}(1,x,x)^{1/2}p^{h}(1,y,y)^{1/2}\quad \forall t>1,\ x,y\in\R .$$
This together with \eqref{1.24} and \eqref{eq:density pht} implies that there is some constant $c_{1}>0$ such that
\begin{equation}\label{1.25}
\left|p^{h}(t,x,y)-1\right|\le c_{1}\e^{-a t}h(x)^{-1}h(y)^{-1}\quad \forall t>1,\ x,y\in\R .
\end{equation}

\begin{lemma}\label{lem7}
Suppose $\mu\in\me$ and $\sigma>0$. For any $f\in\mathcal{B}^{+}_{b}(\R )$ such that $f/h$ is bounded from above and that $\int_{-\infty}^{+\infty }f(x)h(x)dx>0$, we have
$$\lim_{n\to+\infty}\frac{\log \langle f,X_{n\sigma}\rangle}{n\sigma}=\lambda_{1}\quad \pp_{\mu}\mbox{-a.s. on }\{W^{h}_{\infty}(X)>0\}.$$
\end{lemma}

\proof Without loss of generality we assume $0\not=\mu\in\mathcal{M}(\R )$. It follows by Proposition \ref{prop1} that
$$\limsup_{n\to+\infty}\frac{\log \langle f,X_{n\sigma}\rangle}{n\sigma}\le \limsup_{n\to+\infty}\frac{\log\|f\|_{\infty}+\log \mathcal{X}^{0}_{n\sigma}}{n\sigma}\le \lambda_{1}
\quad \pp_{\mu}\mbox{-a.s.}$$
Hence we only need to show that
$$\liminf_{n\to+\infty}\frac{\log \langle f,X_{n\sigma}\rangle}{n\sigma}\ge \lambda_{1}\quad \pp_{\mu}\mbox{-a.s. on }\{W^{h}_{\infty}(X)>0\},$$
or equivalently, for any $\epsilon>0$,
\begin{equation}\label{lem7.2}
\pp_{\mu}\left(\e^{-\lambda_{1}n\sigma}\langle f,X_{n\sigma}\rangle \ge \e^{-\epsilon n\sigma}\text{ for $n$ sufficiently large }|\ W^{h}_{\infty}(X)>0\right)=1.
\end{equation}
For any $n\in\mathbb{N}$ and $\sigma>0$, we have
\begin{align*}
\e^{-\lambda_{1}n\sigma}\langle f,X_{n\sigma}\rangle
&=\mathrm{I}( n,\sigma)+\mathrm{II}( n,\sigma)+\mathrm{III}(n,\sigma),
\end{align*}
where
\begin{align*}
&\mathrm{I}( n,\sigma)=\e^{-\lambda_{1}n\sigma}\langle f,X_{n\sigma}\rangle-\pp_{\mu}\left[\e^{-\lambda_{1}n\sigma}\langle f,X_{n\sigma}\rangle|\mathcal{F}_{n\sigma/2}\right],\\
&\mathrm{II}( n,\sigma)=\pp_{\mu}\left[\e^{-\lambda_{1}n\sigma}\langle f,X_{n\sigma}\rangle|\mathcal{F}_{n\sigma/2}\right]
-(f,h) W^{h}_{n\sigma/2}(X),\\
&\mathrm{III}(n,\sigma)=(f,h) W^{h}_{n\sigma/2}(X).
\end{align*}
Since $\lim_{t\to+\infty}W^{h}_{t}(X)=W^{h}_{\infty}(X)$ $\pp_{\mu}$-a.s., we have
\begin{equation}\label{lem7.3}
\pp_{\mu}\left(\mathrm{III}(n,\sigma)\ge \frac{1}{2}(f,h) W^{h}_{\infty}(X)>0\text{ for $n$ sufficiently large }|\ W^{h}_{\infty}(X)>0\right)=1.
\end{equation}
 Let $\phi(x):=f(x)/h(x)$ for $x\in\R $. By Markov property and \eqref{eq:pht} we have
\begin{align*}
\mathrm{II}( n,\sigma)&=\e^{-\lambda_{1}n\sigma}\langle P^{\beta}_{n\sigma/2}(\phi h),X_{n\sigma/2}\rangle-(\phi h,h) \e^{-\frac{1}{2}\lambda_{1}n\sigma}\langle h,X_{n\sigma/2}\rangle\\
&=\e^{-\frac{1}{2}\lambda_{1}n\sigma}\langle hP^{h}_{n\sigma/2}(\phi),X_{n\sigma/2}\rangle-(\phi h,h) \e^{-\frac{1}{2}\lambda_{1}n\sigma}\langle h,X_{n\sigma/2}\rangle\\
&=\e^{-\frac{1}{2}\lambda_{1}n\sigma}\Big\langle h\int_{-\infty}^{+\infty }\left(p^{h}(n\sigma/2,\cdot,y)-1\right)\phi(y)h^{2}(y)dy,X_{n\sigma/2}\Big\rangle.
\end{align*}
It follows by \eqref{1.25} that for $n\in\mathbb{N}$ with $n\sigma>1$,
$$|\mathrm{II}( n,\sigma)|\le \e^{-\frac{1}{2}\lambda_{1}n\sigma}\langle h\int_{-\infty}^{+\infty }\left|p^{h}(n\sigma/2,\cdot,y)-1\right|\phi(y)h^{2}(y)dy,X_{n\sigma/2}\rangle\le c_{1}\e^{-a n\sigma/2}(\phi,h) \e^{-\frac{1}{2}\lambda_{1}n\sigma} \mathcal{X}^{0}_{n\sigma/2}.$$
This together with \eqref{prop1.2} yields that
\begin{equation}\label{lem7.1}
\pp_{\mu}\left(\lim_{n\to+\infty}|\mathrm{II}( n,\sigma)|=0\right)=1.
\end{equation}
By \eqref{eq:martingale representation} we have
$$\e^{-\lambda_{1}n\sigma}\langle f,X_{n\sigma}\rangle =\langle h P^{h}_{n\sigma}\phi,X_{0}\rangle+\int_{0}^{n\sigma}\int_{-\infty}^{+\infty }\e^{-\lambda_{1}s}h(x)P^{h}_{n\sigma-s}\phi(x)M(ds,dx).$$
Here $[0,n\sigma]\ni t\mapsto \int_{0}^{t}\int_{-\infty}^{+\infty }\e^{-\lambda_{1}s}h(x)P^{h}_{n\sigma-s}\phi(x)M(ds,dx)$ is a square-integrable martingale with quadratic variation $t\mapsto \int_{0}^{t}\langle 2\gamma \e^{-2\lambda_{1}s}h^{2}\left(P^{h}_{n\sigma-s}\phi\right)^{2},X_{s}\rangle ds$.
Hence $$\mathrm{I}(n,\sigma)=\int_{n\sigma/2}^{n\sigma}\int_{-\infty}^{+\infty }\e^{-\lambda_{1}s}h(x)P^{h}_{n\sigma-s}\phi(x)M(ds,dx).$$
Moreover, by \eqref{eq:pht} we have
\begin{align*}
\pp_{\mu}\left[\mathrm{I}(n,\sigma)^{2}\right]
&=\pp_{\mu}\left[\int_{n\sigma/2}^{n\sigma}\langle 2\gamma \e^{-2\lambda_{1}s}h^{2}\left(P^{h}_{n\sigma-s}\phi\right)^{2},X_{s}\rangle ds\right]\\
&\le 2\|\gamma\|_{\infty}\|\phi\|_{\infty}\|h\|_{\infty}\int_{n\sigma/2}^{n\sigma}\e^{-2\lambda_{1}s}\pp_{\mu}\left[\langle hP^{h}_{n\sigma-s}\phi,X_{s}\rangle\right]ds\\
&=c_{2}\int_{n\sigma/2}^{n\sigma}\e^{-2\lambda_{1}s}\langle P^{\beta}_{s}\left(hP^{h}_{n\sigma-s}\phi\right),\mu\rangle ds\\
&=c_{2}\int_{n\sigma/2}^{n\sigma}\e^{-\lambda_{1}s}\langle hP^{h}_{n\sigma}\phi,\mu \rangle ds\\
&\le c_{2}\|\phi\|_{\infty}\langle h,\mu\rangle\int_{n\sigma/2}^{n\sigma}\e^{-\lambda_{1}s}ds\\
&=c_{2}\|\phi\|_{\infty}\langle h,\mu\rangle \lambda_{1}^{-1}\e^{-\lambda_{1}n\sigma/2}\left(1-\e^{-\lambda_{1}n\sigma/2}\right).
\end{align*}
Immediately $\sum_{n=0}^{+\infty}\pp_{\mu}\left[\mathrm{I}(n,\sigma)^{2}\right]<+\infty$. Hence by Borel-Cantelli lemma,
$$\pp_{\mu}\left(\lim_{n\to+\infty}|\mathrm{I}(n,\sigma)|=0\right)=1.$$
This together with \eqref{lem7.3} and \eqref{lem7.1} yields \eqref{lem7.2}. Hence we complete the proof.\qed

\bigskip

\begin{lemma}\label{lem9}
Suppose $0<\delta<\sqrt{\lambda_{1}/2}$
and $\sigma>0$. For any nontrivial $\mu\in\mc$,
$$\liminf_{n\to+\infty}\frac{\log \mathcal{X}^{\delta n\sigma}_{n\sigma}}{n\sigma}\ge -\Lambda_{\delta}\quad \pp_{\mu}\mbox{-a.s. on }\{W^{h}_{\infty}(X)>0\}.$$
\end{lemma}

\proof
We define a quadratic branching mechanism $\tilde{\psi}$ by
$$\tilde{\psi}(x,\lambda):=-\beta(x)\lambda+\gamma(x)\lambda^{2},\quad\forall x\in\R, \lambda\ge 0,$$
where $\gamma$ is defined in \eqref{def-gamma}.
Let $((\tilde{X}_{t})_{t\ge 0},\pp_{\delta_{x}})$ be a $(B_{t},\tilde{\psi})$-superprocess started from Dirac measure at $x$.
For any $R,t\ge 0$ and $x\in\R$, let $\tilde{u}^{R}(t,x):=-\log\pp_{\delta_{x}}\left[\e^{-\tilde{X}((-R,R)^{c})}\right]$ and $u^{R}(t,x):=-\log\pp_{\delta_{x}}\left[\e^{-\mathcal{X}^{R}_{t}}\right]$.
Noting that $\psi\le \tilde{\psi}$, we have by \cite[Corollary 5.18]{Li} that
\begin{equation}\label{lem9.9}
u^{R}(t,x)\ge \tilde{u}^{R}(t,x)\quad\forall t,R\ge 0,x\in\R.
\end{equation}
It is known that $(t,x)\mapsto \tilde{u}^{R}(t,x)$ is the unique nonnegative locally bounded solution to the following integral equation.
\begin{align*}
\tilde{u}^{R}(t,x)&=\Bp_{x}\left(|B_{t}|\ge R\right)+\Bp_{x}\left[\int_{0}^{t}\beta(B_{s})\tilde{u}^{R}(t-s,B_{s})-\gamma(B_{s})\tilde{u}^{R}(t-s,B_{s})^{2}ds\right].
\end{align*}
By \cite[Proposition 2.9]{Li}, $\tilde{u}^{R}(t,x)$ also satisfies that
\begin{equation}\label{lem9.0}
\tilde{u}^{R}(t,x)=\Bp_{x}\left[\exp\left\{\int_{0}^{t}\left(\beta(B_{s})-\gamma(B_{s})\tilde{u}^{R}(t-s,B_{s})\right)ds\right\}\1_{\{|B_{t}|\ge R\}}\right].
\end{equation}
Immediately we have
\begin{equation}\label{lem9.5}
\tilde{u}^{R}(t,x)\le \Bp_{x}\left[e_{\beta}(t),|B_{t}|\ge R\right]=\pi^{R}_{t}(x)\quad \forall x\in\R,\ t\ge 0.
\end{equation}
Let $q\in (0,1)$ and $p=1-q$. By \eqref{lem9.0} and \eqref{lem9.5} one has for all $n\in\mathbb{N}$ and $x\in\R$,
\begin{align}
&\tilde{u}^{\delta n\sigma}(nq\sigma,x)\notag\\
&= \Bp_{x}\left[\exp\left\{\int_{0}^{nq\sigma}\left(\beta(B_{s})-\gamma(B_{s})\tilde{u}^{\delta n\sigma}(nq\sigma-s,B_{s})\right)ds\right\}\1_{\{|B_{nq\sigma}|\ge \delta n\sigma\}}\right]\notag\\
&\ge \Bp_{x}\left[\exp\left\{\int_{0}^{nq\sigma}\left(\beta(B_{s})-\gamma(B_{s})\pi^{\delta n\sigma}_{nq\sigma-s}(B_{s})\right)ds\right\}\1_{\{|B_{nq\sigma}|\ge \delta n\sigma\}}\right]\notag\\
&\ge \Bp_{x}\left[\exp\left\{\int_{0}^{nq\sigma}\left(\beta(B_{s})-\gamma(B_{s})\pi^{\frac{\delta}{q}\left(nq\sigma-s\right)}_{nq\sigma-s}(B_{s})\right)ds\right\}\1_{\{|B_{nq\sigma}|\ge \delta n\sigma\}}\right]\notag\\
&=\Bp_{x}\left[\exp\left\{-\int_{0}^{nq\sigma}\gamma(B_{nq\sigma-r})\pi^{\frac{\delta}{q}r}_{r}(B_{nq\sigma-r})dr\right\}e_{\beta}(nq\sigma)\1_{\{|B_{nq\sigma}|\ge \delta n\sigma\}}\right].\label{lem9.1}
\end{align}
The final equality follows from the changes of variables.
By choosing $q\in (0,\delta/\sqrt{2\lambda_{1}})$, we have $\delta/q>\sqrt{2\lambda_{1}}$.
Since $\gamma$ is compactly supported, it follows from Lemma \ref{lem2}(iii) that
$$c_{1}:=\sup_{x\in\R^{d}}\gamma(x)\int_{0}^{+\infty}\pi^{\frac{\delta}{q}s}_{s}(x)ds<+\infty.$$
Putting this back to \eqref{lem9.1}, one gets
\begin{equation}\label{lem9.2}
\tilde{u}^{\delta n\sigma}(nq\sigma,x)
\ge \e^{-c_{1}}\Bp_{x}\left[e_{\beta}(nq\sigma)\1_{\{|B_{nq\sigma}|\ge \delta n\sigma\}}\right]
=\e^{-c_{1}}\pi^{\frac{\delta}{q}\cdot nq\sigma}_{nq\sigma}(x).
\end{equation}
Let $K$ be a compact set of $\R$ with $\int_{K}h(x)dy>0$.
Since $\frac{\delta}{q}>\sqrt{2\lambda_{1}}$, we have by Lemma \ref{lem1} that
$$\lim_{n\to+\infty}\inf_{x\in K}\frac{\log \pi^{\delta n\sigma}_{nq\sigma}(x)}{nq\sigma}=-\Lambda_{\frac{\delta}{q}}=-\frac{\delta^{2}}{2q^{2}}.$$
Let $\epsilon\in (0,-\Lambda_{\delta}/8)$.
The above equation combined with \eqref{lem9.9} and \eqref{lem9.2} implies that
there is some constant $c_{2}>0$ such that
\begin{equation}\label{lem9.6}
\inf_{x\in K}u^{\delta n\sigma}(nq\sigma,x)\ge c_{2}\e^{-\frac{\delta^{2}}{2q}n\sigma- \frac{1}{2}\epsilon nq\sigma}
\end{equation}
for $n$ sufficiently large.
By the Markov property, we have
\begin{align*}
&\pp_{\mu}\left(\mathcal{X}^{\delta n\sigma}_{n\sigma}\le \e^{\left(-\Lambda_{\delta}- \frac{3}{2}\epsilon\right)n\sigma};\ X_{np\sigma}(K)\ge\e^{\left(\lambda_{1}- \frac{1}{2}\epsilon\right)np\sigma}\right)\\
&=\pp_{\mu}\left[\pp_{X_{np\sigma}}\left[\e^{-\mathcal{X}^{\delta n\sigma}_{nq\sigma}}\ge \e^{-\e^{(-\Lambda_{\delta}-\frac{3}{2}\epsilon)n\sigma}}\right];X_{np\sigma}(K)\ge \e^{(\lambda_{1}-\frac{1}{2}\epsilon)np\sigma}\right]\\
&\le \exp\{\e^{\left(-\Lambda_{\delta}- \frac{3}{2}\epsilon\right)n\sigma}\}\pp_{\mu}\left[\e^{-\langle u^{\delta n\sigma}(nq\sigma,\cdot),X_{np\sigma}\rangle};\ X_{np\sigma}(K)\ge\e^{\left(\lambda_{1}- \frac{1}{2}\epsilon\right)np\sigma}\right]\\
&\le \exp\{\e^{\left(-\Lambda_{\delta}- \frac{3}{2}\epsilon\right)n\sigma}\}\pp_{\mu}\left[\exp\{-c_{2}\e^{-\frac{\delta^{2}}{2q}n\sigma- \frac{1}{2}\epsilon nq\sigma}X_{n p\sigma}(K)\};\ X_{np\sigma}(K)\ge\e^{\left(\lambda_{1}- \frac{1}{2}\epsilon\right)np\sigma}\right]\\
&\le \exp\{\e^{\left(-\Lambda_{\delta}- \frac{3}{2}\epsilon\right)n\sigma}-c_{2}\e^{-\frac{\delta^{2}}{2q}n\sigma-\frac{1}{2}\epsilon nq\sigma+(\lambda_{1}-\frac{1}{2}\epsilon)np\sigma}\}
\pp_{\mu}\left[X_{np\sigma}(K)\ge\e^{\left(\lambda_{1}- \frac{1}{2}\epsilon\right)np\sigma}\right]\\
&= \exp\left\{-\e^{\left(-\Lambda_{\delta}- \frac{3}{2}\epsilon\right)n\sigma}\left(c_{2}\e^{n\sigma\left(\sqrt{2\lambda_{1}}\delta+\epsilon-\lambda_{1}q-\frac{\delta^{2}}{2q}\right)}-1\right)\right\}\pp_{\mu}\left[X_{np\sigma}(K)\ge\e^{\left(\lambda_{1}- \frac{1}{2}\epsilon\right)np\sigma}\right].
\end{align*}
The first and second inequalities follow from Chebyshev's inequality and \eqref{lem9.6}, respectively.
We choose
$q\in \left(\left(\frac{\delta}{\sqrt{2\lambda_{1}}}+\frac{\epsilon-\sqrt{\epsilon^{2}+2\sqrt{2\lambda_{1}}\delta\epsilon}}{2\lambda_{1}}\right)\vee 0,\frac{\delta}{\sqrt{2\lambda_{1}}}\right)$
then $\sqrt{2\lambda_{1}}\delta+\epsilon-\lambda_{1}q-\frac{\delta^{2}}{2q}>0$. Note that $-\Lambda_{\delta}- \frac{3}{2}\epsilon>-13\Lambda_{\delta}/16>0$ for $0<\delta<\sqrt{\lambda_{1}/2}$ and $\epsilon\in (0,-\Lambda_{\delta}/8)$. The above inequality implies that $\pp_{\mu}\left(\mathcal{X}^{\delta n\sigma}_{n\sigma}\le \e^{\left(-\Lambda_{\delta}- \frac{3}{2}\epsilon\right)n\sigma};\ X_{np\sigma}(K)\ge\e^{\left(\lambda_{1}- \frac{1}{2}\epsilon\right)np\sigma}\right)$ decreases faster than exponentially as $n\to+\infty$. Thus
$$\sum_{n=0}^{+\infty}\pp_{\mu}\left(\mathcal{X}^{\delta n\sigma}_{n\sigma}\le \e^{\left(-\Lambda_{\delta}- \frac{1}{2}\epsilon\right)n\sigma};\ X_{np\sigma}(K)\ge\e^{\left(\lambda_{1}- \frac{1}{2}\epsilon\right)np\sigma}\right)<+\infty,$$
and by Borel-Cantelli lemma,
\begin{equation}\label{lem9.3}
\pp_{\mu}\left(\frac{\log \mathcal{X}^{\delta n\sigma}_{n\sigma}}{n\sigma}>-\Lambda_{\delta}- \frac{1}{2}\epsilon\mbox{ or }\frac{\log X_{np\sigma}(K)}{np\sigma}<\lambda_{1}- \frac{1}{2}\epsilon\quad\mbox{ for $n$ sufficiently large}\right)=1.
\end{equation}
Note that by Lemma \ref{lem7}
\begin{equation}\nonumber
\pp_{\mu}\left(\lim_{n\to+\infty}\frac{\log X_{np\sigma}(K)}{np\sigma}=\lambda_{1}\Big|W^{h}_{\infty}(X)>0\right)=1.
\end{equation}
We get by \eqref{lem9.3} that
$$\pp_{\mu}\left(\frac{\log \mathcal{X}^{\delta n\sigma}_{n\sigma}}{n\sigma}>-\Lambda_{\delta}- \frac{1}{2}\epsilon\quad\mbox{ for $n$ sufficiently large }\Big| W^{h}_{\infty}(X)>0\right)=1.$$
This lemma follows by letting $\epsilon\downarrow 0$.\qed

\bigskip

\begin{proposition}\label{prop3}
For $0\le \delta<\sqrt{\lambda_{1}/2}$ and $\mu\in\mc$,
$$\liminf_{t\to+\infty}\frac{\log \mathcal{X}^{\delta t}_{t}}{t}\ge -\Lambda_{\delta}\quad \pp_{\mu}\mbox{-a.s. on }\{W^{h}_{\infty}(X)>0\}.$$
\end{proposition}

\proof  First we consider $\delta=0$. Since
$$W^{h}_{t}(X)=\e^{-\lambda_{1}t}\langle h,X_{t}\rangle \le \|h\|_{\infty}\e^{-\lambda_{1}t}\mathcal{X}^{0}_{t},$$
we have
$$\frac{\log \mathcal{X}^{0}_{t}}{t}\ge \frac{\log W^{h}_{t}(X)}{t}-\frac{\log \|h\|_{\infty}}{t}+\lambda_{1}.$$
Since $W^{h}_{t}(X)\to W^{h}_{\infty}(X)$ $\pp_{\mu}$-a.s., we get that
$$\liminf_{t\to+\infty}\frac{\log \mathcal{X}^{0}_{t}}{t}\ge \lambda_{1}=-\Lambda_{0}\quad\pp_{\mu}\mbox{-a.s. on }\{W^{h}_{\infty}(X)>0\}.$$

Now suppose $0<\delta<\sqrt{\lambda_{1}/2}$.
Let $\sigma>0$. We take $\theta>0$ small such that $\delta_{\theta}:=\delta+\theta<\sqrt{\lambda_{1}/2}.$
By \eqref{fact:6} we have
$$\pp_{\mu}\left[\mathcal{X}^{\delta_{\theta}(n+1)\sigma}_{(n+1)\sigma}\Big|\mathcal{F}_{n\sigma}\right]
\ge c_{1}\mathcal{X}^{\delta_{\theta}n\sigma}_{n\sigma}$$
for some constant $c_{1}=c_{1}(\theta,\sigma)>0.$ It then follows by Lemma \ref{lem9} that
\begin{equation}\label{lem10.1}
\pp_{\mu}\left(\liminf_{n\to+\infty}\frac{\log \pp_{\mu}\left[\mathcal{X}^{\delta_{\theta}(n+1)\sigma}_{(n+1)\sigma}|\mathcal{F}_{n\sigma}\right]}{n\sigma}\ge -\Lambda_{\delta_{\theta}}\Big|W^{h}_{\infty}(X)>0\right)=1.
\end{equation}
Let $0<\epsilon<-\Lambda_{\delta_{\theta}}/4$.
We have
\begin{align}
&\sup_{t\in [n\sigma,(n+1)\sigma)}e^{\Lambda_{\delta_{\theta}}t}\left|\pp_{\mu}\left[\mathcal{X}^{\delta_{\theta}(n+1)\sigma}_{(n+1)\sigma}|\mathcal{F}_{t}\right]
-\pp_{\mu}\left[\mathcal{X}^{\delta_{\theta}(n+1)\sigma}_{(n+1)\sigma}|\mathcal{F}_{n\sigma}\right]\right|\notag\\
&\le \sup_{t\in [n\sigma,(n+1)\sigma)}e^{\Lambda_{\delta_{\theta}}n\sigma}\left|\pp_{\mu}\left[\mathcal{X}^{\delta_{\theta}(n+1)\sigma}_{(n+1)\sigma}|\mathcal{F}_{t}\right]
-\pp_{\mu}\left[\mathcal{X}^{\delta_{\theta}(n+1)\sigma}_{(n+1)\sigma}|\mathcal{F}_{n\sigma}\right]\right|\notag\\
&= e^{-\left(-\frac{1}{2}\Lambda_{\delta_{\theta}}-\epsilon\right)n\sigma-\left(\frac{1}{2}\Lambda_{\delta_{\theta}}-\epsilon\right)\sigma} \sup_{t\in [n\sigma,(n+1)\sigma)}e^{\left(\frac{1}{2}\Lambda_{\delta_{\theta}}-\epsilon\right)(n+1)\sigma}\left|\pp_{\mu}\left[\mathcal{X}^{\delta_{\theta}(n+1)\sigma}_{(n+1)\sigma}|\mathcal{F}_{t}\right]
-\pp_{\mu}\left[\mathcal{X}^{\delta_{\theta}(n+1)\sigma}_{(n+1)\sigma}|\mathcal{F}_{n\sigma}\right]\right|.\notag
\end{align}
By Lemma \ref{lem5}, the final term in the right hand side converges to $0$ as $n\to+\infty$. Thus we get
\begin{equation}\label{lem10.2}
\lim_{n\to+\infty}\sup_{t\in [n\sigma,(n+1)\sigma)}e^{\Lambda_{\delta_{\theta}}t}\left|\pp_{\mu}\left[\mathcal{X}^{\delta_{\theta}(n+1)\sigma}_{(n+1)\sigma}|\mathcal{F}_{t}\right]
-\pp_{\mu}\left[\mathcal{X}^{\delta_{\theta}(n+1)\sigma}_{(n+1)\sigma}|\mathcal{F}_{n\sigma}\right]\right|=0\quad\pp_{\mu}\mbox{-a.s.}
\end{equation}
Note that for any $t\in [n\sigma,(n+1)\sigma)$,
\begin{align*}
\e^{\Lambda_{\delta_{\theta}}t}\pp_{\mu}\left[\mathcal{X}^{\delta_{\theta}(n+1)\sigma}_{(n+1)\sigma}|\mathcal{F}_{t}\right]
&\ge \e^{\Lambda_{\delta_{\theta}}(n+1)\sigma}\pp_{\mu}\left[\mathcal{X}^{\delta_{\theta}(n+1)\sigma}_{(n+1)\sigma}|\mathcal{F}_{n\sigma}\right]\\
&\quad -\sup_{t\in [n\sigma,(n+1)\sigma)}e^{\Lambda_{\delta_{\theta}}t}\left|\pp_{\mu}\left[\mathcal{X}^{\delta_{\theta}(n+1)\sigma}_{(n+1)\sigma}|\mathcal{F}_{t}\right]
-\pp_{\mu}\left[\mathcal{X}^{\delta_{\theta}(n+1)\sigma}_{(n+1)\sigma}|\mathcal{F}_{n\sigma}\right]\right|
\end{align*}
Hence by \eqref{lem10.1} and \eqref{lem10.2} we get that
\begin{equation}\label{lem10.3}
\pp_{\mu}\left(\liminf_{t\to+\infty}\frac{\log \pp_{\mu}\left[\mathcal{X}^{\delta_{\theta}(n+1)\sigma}_{(n+1)\sigma}|\mathcal{F}_{t}\right]}{t}\ge -\Lambda_{\delta_{\theta}}\Big|W^{h}_{\infty}(X)>0\right)=1.
\end{equation}
By Markov property, for any $t\in [n\sigma,(n+1)\sigma)$,
\begin{align}
\pp_{\mu}\left[\mathcal{X}^{\delta_{\theta}(n+1)\sigma}_{(n+1)\sigma}|\mathcal{F}_{t}\right]
&=\langle \pi^{\delta_{\theta}(n+1)\sigma}_{(n+1)\sigma-t},X_{t}\rangle.\notag
\end{align}
So we have
\begin{align}
\mathrm{I}(\delta_{\theta},t):=\langle \pi^{\delta_{\theta}(n+1)\sigma}_{(n+1)\sigma-t}\1_{(-\delta t,\delta t)^{c}},X_{t}\rangle
&=\pp_{\mu}\left[\mathcal{X}^{\delta_{\theta}(n+1)\sigma}_{(n+1)\sigma}|\mathcal{F}_{t}\right]-\langle \pi^{\delta_{\theta}(n+1)\sigma}_{(n+1)\sigma-t}\1_{(-\delta t,\delta t)},X_{t}\rangle\notag\\
&=:\pp_{\mu}\left[\mathcal{X}^{\delta_{\theta}(n+1)\sigma}_{(n+1)\sigma}|\mathcal{F}_{t}\right]-\mathrm{II}(\delta_{\theta},t).\label{lem10.5}
\end{align}
Lemma \ref{lem2}(i) implies that there is a constant $c_{2}>0$ independent of $\delta_{\theta}$ and $\theta$ such that
$$\pi^{\delta_{\theta}(n+1)\sigma}_{(n+1)\sigma-t}(y)\le c_{2}(\theta(n+1)\sigma)^{-1}\e^{-\frac{\theta^{2}}{2}(n+1)^{2}\sigma}\quad \forall t\in [n\sigma,(n+1)\sigma),\ |y|<\delta t,
$$
when $n$ is sufficiently large. Hence we get that for $t\in [n\sigma,(n+1)\sigma)$,
\begin{align*}
\mathrm{II}(\delta_{\theta},t)&\le c_{2}(\theta(n+1)\sigma)^{-1}\e^{-\frac{\theta^{2}}{2}(n+1)^{2}\sigma}\langle \1_{(-\delta t,\delta t)},X_{t}\rangle\\
&\le \frac{c_{2}(\theta(n+1)\sigma)^{-1}\e^{-\frac{\theta^{2}}{2}(n+1)^{2}\sigma+\lambda_{1}t}}{\inf_{|y|<\delta(n+1)\sigma}h(y)}
\langle \e^{-\lambda_{1}t}h\1_{(-\delta t,\delta t)},X_{t}\rangle\\
&\le \frac{c_{2}(\theta(n+1)\sigma)^{-1}\e^{-\frac{\theta^{2}}{2}(n+1)^{2}\sigma+\lambda_{1}t}}{\inf_{|y|<\delta(n+1)\sigma}h(y)}W^{h}_{t}(X).
\end{align*}
By \eqref{estimateforh}, there is a constant $c_{3}>0$ such that when $n$ is sufficiently large,
\begin{equation*}
\mathrm{II}(\delta_{\theta},t)\le c_{3}(\theta(n+1)\sigma)^{-1}\e^{-\frac{\theta^{2}}{2}(n+1)^{2}\sigma+(\lambda_{1}+\sqrt{2\lambda_{1}}\delta)(n+1)\sigma}
W^{h}_{t}(x)\quad\forall t\in [n\sigma,(n+1)\sigma)
\end{equation*}
This implies that
\begin{equation}\label{lem10.4}
\lim_{n\to+\infty}\sup_{t\in [n\sigma,(n+1)\sigma)}\e^{\Lambda_{\delta_{\theta}t}}\mathrm{II}(\delta_{\theta},t)=0\quad \pp_{\mu}\mbox{-a.s.}
\end{equation}
Note that by \eqref{lem10.5}
\begin{align*}
\e^{\Lambda_{\delta_{\theta}}t}\mathrm{I}(\delta_{\theta},t)&=\e^{t\left(\log \mathrm{I}(\delta_{\theta},t)/t+\Lambda_{\delta_{\theta}}\right)}\\
&= \e^{t\left(\log \pp_{\mu}\left[\mathcal{X}^{\delta_{\theta}(n+1)\sigma}_{(n+1)\sigma}\big|\mathcal{F}_{t}\right]\big/t+\Lambda_{\delta_{\theta}}\right)}
-\e^{\Lambda_{\delta_{\theta}}t}\mathrm{II}(\delta_{\theta},t).
\end{align*}
This together with \eqref{lem10.3} and \eqref{lem10.4} implies that
\begin{equation}\label{limitinf-I}
\pp_{\mu}\left(\liminf_{t\to+\infty}\frac{\log \mathrm{I}(\delta_{\theta},t)}{t}\ge -\Lambda_{\delta_{\theta}} \Big| W^{h}_{\infty}(X)>0\right)=1.
\end{equation}
Note that by definition
$$\mathrm{I}(\delta_{\theta},t)\le \|\pi^{\delta_{\theta}(n+1)\sigma}_{(n+1)\sigma-t}\|_{\infty}\mathcal{X}^{\delta t}_{t}\le \e^{\|\beta^{+}\|_{\infty}\sigma}\mathcal{X}^{\delta t}_{t}\quad\forall t\in [n\sigma,(n+1)\sigma).$$
By \eqref{limitinf-I}, we have
$$\pp_{\mu}\left(\liminf_{t\to+\infty}\frac{\log \mathcal{X}^{\delta t}_{t}}{t}\ge -\Lambda_{\delta_{\theta}}| W^{h}_{\infty}(X)>0\right)=1.$$\qed

\bigskip

\noindent\textit{Proof of Theorem \ref{them1}:} Theorem \ref{them1} follows immediately from Propositions \ref{prop1} and \ref{prop3}.\qed

\section{Proofs of Theorem \ref{them2} and Theorem \ref{them3}}\label{sec4}

\subsection{Skeleton decomposition}
 It is well-known that the $(B_{t},\psi)$-superprocess can be constructed as the high density limit of a sequence of branching Markov processes. Another link between superprocesses and branching Markov processes is provided by the so-called skeleton decomposition, which is developed by
 \cite{CRY2,EKW,KPR}.
 The skeleton decomposition provides a pathwise description of
a superprocesses in terms of immigrations along a branching Markov process called the skeleton.
The following condition is fundamental for the skeleton construction.

\medskip

There is a locally bounded function $w>0$ on $\R$ satisfying that
\begin{equation}\label{a2}
\pp_{\mu}\left[\e^{-\langle w,X_{t}\rangle}\right]=\e^{-\langle w,\mu\rangle}\quad\forall \mu\in \mc.
\tag{A2}
\end{equation}
This locally bounded martingale function $w$ assures that
$\left(w(B_{t})\exp\{-\int_{0}^{t}\frac{\psi(B_{s},w(B_{s}))}{w(B_{s})}ds\}\right)_{t\ge 0}$ is a $\Bp_{x}$-(super)martingale. Thus one can define a family of (sub)probability measures $\{\Bp^{w}_{x},x\in\R\}$
by
\begin{equation}
\left.\frac{d\Bp^{w}_{x}}{d\Bp_{x}}\right|_{\sigma(B_{s}:s\in [0,t])}:=\frac{w(B_{t})}{w(x)}\exp\left\{-\int_{0}^{t}\frac{\psi(B_{s},w(B_{s}))}{w(B_{s})}ds\right\}\quad\forall t\ge 0.
\end{equation}
We denote the process $((B_{t})_{t\ge 0},\Bp^{w}_{x},x\in\R)$ by $(B^{w}_{t})_{t\ge 0}$.

An integer-valued locally finite random measure $\xi$ on $\R$ is called a point process. If there is a locally finite measure $\lambda$ on $\R$ such that
$\xi(B)$ is
Poisson distributed
with mean $\lambda(B)$ for any Borel set $B$, and that $\xi(B_1),\cdots, \xi(B_n)$ are independent
for any disjoint Borel sets $B_1,\cdots, B_n$, $n\geq 2$,
then $\xi$ is called \textit{Poisson point process} with intensity $\lambda$.
If we randomize by replacing the fixed measure $\lambda$ by
a random measure $\Lambda$ on $\R$,  then we get a \textit{Cox process} directed by $\Lambda$. More precisely,
given $\Lambda$, $\xi$ is conditionally Poisson with intensity $\Lambda$ almost surely.

\begin{proposition}\label{propskeleton}
Assume \eqref{a2} holds.
For every $\mu\in\mf$ there exists a probability space with probability measure $\p_{\mu}$
that carries two processes $(Z_{t})_{t\ge 0}$ and
 $(\widehat{X}_{t})_{t\ge 0}$
satisfying the following conditions.
\begin{description}
\item{(i)} $((Z_{t})_{t\ge 0},\p_{\mu})$ is branching Markov process
with $Z_{0}$ being a Poisson point process with intensity $w(x)\mu(dx)$,
in which
each particle moves independently as a copy of $(B^{w}_{t})_{t\ge 0}$,
and a particle at location $x$ dies at rate $q(x)$ and is replaced by a random number of offspring with distribution $\{p_{k}(x):k\ge 2\}$ uniquely identified by
\begin{eqnarray*}
G(x,s)&:=&q(x)\sum_{k=2}^{+\infty}p_{k}(x)(s^{k}-s)\\
&=&\frac{1}{w(x)}\left[\psi(x,w(x)(1-s))-(1-s)\psi(x,w(x))\right].
\end{eqnarray*}
\item{(ii)}
$((\widehat{X}_{t})_{t\ge 0},\p_{\mu})$ has the same  distribution  as $(X,\pp_{\mu})$.
\item{(iii)}
For every $t\ge 0$,
$Z_{t}$ is a Cox process directed by
 $w\widehat{X}_{t}$.
\end{description}
\end{proposition}

We show in the next proposition that
the martingale function $w$ in \eqref{a2}
exists for the $(B_{t},\psi)$-superprocess.

\begin{proposition}
Let $\mathcal{E}:=\{ W^{h}_{\infty}(X)=0\}$ and $w(x):=-\log\pp_{\delta_{x}}\left(\mathcal{E}\right)$ for $x\in\R$. Then $w$ is a
bounded positive function
satisfying \eqref{a2}.
Moreover $w'(x)=0$ for $|x|$ sufficiently large.
\end{proposition}

\proof
Since $ W^{h}_{\infty}(X)$ is nondegenerate under $\pp_{\delta_{x}}$,
$w(x)=-\log\pp_{\delta_{x}}\left(\mathcal{E}\right)$ takes values in $(0,+\infty]$.
 We only need to show that $w$ is a bounded function on $\R$ since the second assertion is a direct result of Lemma \ref{lem4.1} and the boundedness of $w$.

Let $\mathcal{E}_{ext}:=\{\|X_{t}\|=0\mbox{ for $t$ sufficiently large}\}$ and $w_{ext}(x):=-\log\pp_{\delta_{x}}\left(\mathcal{E}_{ext}\right)$ for $x\in\R$.
Since $\psi(x,\lambda)\ge -\beta(x)\lambda+\alpha(x)\lambda^{2}=:\hat{\psi}(x,\lambda)$ for $x\in\R$ and $\lambda\ge 0$, it follows by \cite[Corollary 5.18]{Li} that the extinction probability of the $(B_{t},\psi)$-superprocess is larger than that of the $(B_{t},\hat{\psi})$-superprocess.
Let $\hat{w}_{ext}$ be the log-Laplace exponent of the $(B_{t},\hat{\psi})$-superprocess. Let $\emptyset\not=\mathcal{O}:=\{x\in\R:\alpha(x)>0\}$.
By \cite[Lemma 7.1]{EP99}, $\hat{w}_{ext}(x)$ is a locally bounded function on $\mathcal{O}$. Since $\mathcal{E}_{ext}\subseteq\mathcal{E}$, one has $w(x)\le w_{ext}(x)\le\hat{w}_{ext}(x)<+\infty$ for all $x\in\mathcal{O}$.
In the remaining of this proof, we fix an arbitrary $c\in\mathcal{O}$ and $x_{0}\in \R\setminus\mathcal{O}$. Without loss of generality, we assume $x_{0}>c$.

For $t\in (0,+\infty)$, let $Q_{t}:=(0,t)\times (c,+\infty)$ be the time-space open set. Let $X^{Q_{t}}$ be the exit measure from $Q_{t}$. Then for any initial measure $\mu$ with support contained in $(c,+\infty)$, $X^{Q_{t}}$ is a finite random measure supported on the boundary of $Q_{t}$. Let $\mathcal{H}$ be the set of nonnegative bounded functions on $[0,+\infty)\times\R$ satisfying that there is some $S$ such that $f(s,y)=0$ for all $(s,y)\in [S,+\infty)\times\R$. For $f\in\mathcal{H}$, let $U_{f}(t,x):=-\log\pp_{\delta_{x}}\left[\exp\{-\langle f,X^{Q_{t}}\rangle\}\right]$ for $t\ge 0$ and $x\in\R$. Then $U_{f}(t,x)$ satisfies the following integral equation.
\begin{equation}\label{lem4.2.1}
U_{f}(t,x)+\Bp_{x}\left[\int_{0}^{\tau_{c}\wedge t}\psi(B_{s},U_{f}(t-s,B_{s}))\right]=\Bp_{x}\left[f(\tau_{c}\wedge t,B_{\tau_{c}\wedge t})\right],
\end{equation}
where $\tau_{c}$ denotes the first exit time of $(B_{t})_{t\ge 0}$ from $(c,+\infty)$.
Let $X^{c}_{t}(A):=X^{Q_{t}}(\{t\}\times (A\cap (c,+\infty)))$ for any $A\subset\R$. This definition implies that $X^{c}_{t}$ is the projection of $X^{Q_{t}}$ on $\{t\}\times (c,+\infty)$. Let $u^{c}_{g}(t,x):=-\log\pp_{\delta_{x}}\left[\exp\{-\langle g,X^{c}_{t}\rangle\}\right]$ for $g\in\mathcal{B}^{+}_{b}(\R)$, $t\ge 0$ and $x\in \R$. By \eqref{lem4.2.1}, $u^{c}_{g}(t,x)$ satisfies the following integral equation.
\begin{equation}\nonumber
u^{c}_{g}(t,x)+\Bp_{x}\left[\int_{0}^{t}\psi(B^{c}_{s},u^{c}_{g}(t-s,B^{c}_{s}))ds\right]=\Bp_{x}\left[g(B^{c}_{t})\right],
\end{equation}
where $(B^{c}_{t})_{t\ge 0}$ denotes the Brownian motion killed outside $(c,+\infty)$. This implies that $(X^{c}_{t})_{t\ge 0}$ is a $(B^{c}_{t},\psi)$-superprocess.
Note that
\begin{eqnarray}
\pp_{\delta_{x_{0}}}\left[\e^{-\lambda_{1}t}\langle h,X^{c}_{t}\rangle\right]&=&\e^{-\lambda_{1}t}\Bp_{x_{0}}\left[e_{\beta}(t)h(B_{t});t<\tau_{c}\right]\nonumber\\
&=&\Bp^{h}_{x_{0}}\left(t<\tau_{c}\right).
\end{eqnarray}
Here $\Bp^{h}_{x_{0}}$ is the probability measure defined by
$$\left.\frac{d\Bp^{h}_{x_{0}}}{d\Bp_{x_{0}}}\right|_{\sigma(B_{s}:s\le t)}=\e^{-\lambda_{1}t}e_{\beta}(t)\frac{h(B_{t})}{h(x_{0})}\quad\forall t\ge 0.$$
It is known that $((B_{t})_{t\ge 0},\Bp^{h}_{x_{0}})$ is a recurrent diffusion on $\R$. So $\Bp^{h}_{x_{0}}\left(t<\tau_{c}\right)\to 0$ as $t\to+\infty$.
This implies that $\e^{-\lambda_{1}t}\langle h,X^{c}_{t}\rangle$ converges to $0$ in $L^{1}(\pp_{\delta_{x_{0}}})$, and so there is a subsequence of $\e^{-\lambda_{1}t}\langle h,X^{c}_{t}\rangle$ which converges to $0$ $\pp_{\delta_{x_{0}}}$-a.s.

On the other hand, we note that $\|X^{Q_{t}}\left((0,+\infty)\times \{c\}\right)\|$ denotes the total mass of the projection of $X^{Q_{t}}$ on $(0,t]\times \{c\}$. For $\lambda, t\ge 0$ and $y\in\R$, let $v^{c}_{\lambda}(t,y):=-\log\pp_{\delta_{y}}\left[\exp\{-\lambda\|X^{Q_{t}}((0,+\infty)\times\{c\})\|\}\right]$. It follows by \eqref{lem4.2.1} that
\begin{eqnarray}
v^{c}_{\lambda}(t,x)&=&\lambda\Bp_{x}\left(\tau_{c}\le t\right)-\Bp_{x}\left[\int_{0}^{\tau_{c}\wedge t}\psi(B_{s},v^{c}_{\lambda}(t-s,B_{s}))ds\right]\nonumber\\
&=&\lambda\Bp_{x}\left[e_{\beta}(\tau_{c}\wedge t);\tau_{c}\le t\right]-\Bp_{x}\left[\int_{0}^{\tau_{c}\wedge t}e_{\beta}(s)\psi_{0}(B_{s},v^{c}_{\lambda}(t-s,B_{s}))ds\right],\label{lem4.2.2}
\end{eqnarray}
where $\psi_{0}(x,\lambda)=\psi(x,\lambda)+\beta(x)\lambda$.
The second equation follows from \cite[Lemma A.1]{EKW}. Let $\|X^{\{c\}}\|$ be the limit of the nondecreasing sequence $\{\|X^{Q_{t}}((0,+\infty)\times\{c\})\|:\ t\ge 0\}$ and $v^{c}_{\lambda}(x):=-\log\pp_{\delta_{x}}\left[\exp\{-\lambda\|X^{\{c\}}\|\}\right]=\lim_{t\to+\infty}v^{c}_{\lambda}(t,x)$ for $\lambda\ge 0$ and $x\ge c$. By \eqref{lem4.2.2}, one has $v^{c}_{\lambda}(t,x)\le \lambda\Bp_{x}\left[e_{\beta}(\tau_{c});\tau_{c}\le t\right]$, and so $v^{c}_{\lambda}(x)\le \lambda \Bp_{x}\left[e_{\beta}(\tau_{c})\right]$. Since $\beta$ is compactly supported, by \cite[Theorem 9.22]{CZ} $x\mapsto \Bp_{x}\left[e_{\beta}(\tau_{c})\right]$ is a bounded function on $[c,+\infty)$. Thus for every $\lambda\ge 0$, $x\mapsto v^{c}_{\lambda}(x)$ is a bounded function on $[c,+\infty)$.

We note that $\mathcal{E}=\{ W^{h}_{\infty}(X)=0\}=\{\exists\ t_{n}\to+\infty \mbox{ such that } \e^{-\lambda_{1}t_{n}}\langle h,X_{t_{n}}\rangle\to 0\}$. We have
\begin{eqnarray}
\e^{-w(x_{0})}=\pp_{\delta_{x_{0}}}(\mathcal{E})&=&\pp_{\delta_{x_{0}}}\left[\pp_{\delta_{x_{0}}}\left(\mathcal{E}|X^{Q_{t}}\right)\right]\nonumber\\
&\ge&\pp_{\delta_{x_{0}}}\left[\pp_{\|X^{\{c\}}\|\delta_{c}}\left(\mathcal{E}_{ext}\right)\right]\nonumber\\
&=&\pp_{\delta_{x_{0}}}\left[\e^{-w_{ext}(c)\|X^{\{c\}}\|}\right]=\e^{-v^{c}_{w_{ext}(c)}(x_{0})}.\nonumber
\end{eqnarray}
Thus one gets $w(x_{0})\le v^{c}_{w_{ext}(c)}(x_{0})$ and so $w$ is a bounded function on $[c,+\infty)$.\qed

\bigskip

\subsection{Limiting distributions for the skeleton}

Since
$({\widehat X};\p_{\mu})$
is equal in distribution to the $(B_{t},\psi)$-superprocess, we may work on this skeleton space whenever it is convenient.
For notational simplification, we will abuse the notation and denote $\widehat X$ by $X$.
We will refer to $(Z_{t})_{t\ge 0}$ as the skeleton branching diffusion (skeleton) of $X$.
We use $u\in Z_{t}$ to denote a particle of the skeleton which is alive at time t, and $z_{u}(t)$ for its spatial location.
We use $\|Z_{t}\|$ to denote the total number of particles alive at time $t$.

In this section we shall show that the skeleton branching diffusion $Z_{t}$
shifted by $\sqrt{\lambda_{1}/2}\,t$
converges in distribution to a Cox process directed by a random measure which has a random intensity mixed by the limit of
an additive martingale
(see Proposition \ref{prop5} below). Our proof follows the same approach as \cite{Bocharov} (see, also \cite{Nishimori}):
First, we represent
the population moments
in terms of Feynman-Kac functionals associated to Brownian motions,
see \eqref{meanZ} and \eqref{varianceZ} below.
Using the estimates established in Section \ref{sec2}, we show in Lemma \ref{lem4.4} that the second order moment is asymptotically the same as the first order moment. Combining this with the Chebyshev and Payley-Zigmund inequalities, we compute the asymptotic behavior of the distributions of particles near $\sqrt{\lambda_{1}/2}\,t$ in Lemma \ref{lem4.5}. We can then follow the argument of \cite{Bocharov} to establish Proposition \ref{prop4}.

Recall that $w'(x)=0$ when $|x|$ is large. So we assume that there are constants $M,w_{\pm}>0$ such that $w(x)=w_{-}$ for $x\ge M$ and $w(x)=w_{+}$ for $x\le -M$.

In what follows we always assume the following:
\begin{description}
\item{(1)} $R(t)=\delta t+a(t)$ where $\delta\in (0,\sqrt{2\lambda_{1}})$ and $a(t)=o(t)$ as $t\to+\infty$.
\item{(2)} For some $a\in \left(0,1-\frac{\delta}{\sqrt{2\lambda_{1}}}\right)$, $0\le s(t)<a t$ for all $t\ge 0$ and $s(t)=o(t)$ as $t\to+\infty$.
\item{(3)} $b(t)\ge 0$ for all $t\ge 0$ and $b(t)=o(t)$ as $t\to+\infty$.
\item{(4)} $x(\cdot): [0,+\infty)\to\R$ satisfies $|x(t)|\le b(t)$ for all large $t>0$.
\item{(5)} $A$ is a Borel set of $\R$ with $\inf A>-\infty$.
\end{description}

We use $\bp_{\nu}$ to denote the probability measure where the branching Markov process $(Z_{t})_{t\ge 0}$ started from the integer-valued measure $\nu$.
For every $x\in\R$, the fist two moments of $((Z_{t})_{t\ge 0},\bp_{\delta_{x}})$
can be expressed by the spatial motion and the branching rate
: For $f\in\mathcal{B}_{b}^{+}(\R)$,
\begin{eqnarray}
\bp_{\delta_{x}}\left[\langle f,Z_{t}\rangle\right]&=&\Bp^{w}_{x}\left[\e^{\int_{0}^{t}\frac{\partial}{\partial s}G(B^{w}_{r},1)dr}f(B^{w}_{t})\right]=\frac{1}{w(x)}P^{\beta}_{t}(w f)(x).\label{meanZ}\\
\bp_{\delta_{x}}\left[\langle f,Z_{t}\rangle^{2}\right]&=&\bp_{\delta_{x}}\left[\langle f^{2},Z_{t}\rangle\right]+\Bp^{w}_{x}\left[\int_{0}^{t}\e^{\int_{0}^{r}\frac{\partial}{\partial s}G(B^{w}_{u},1)du}\frac{\partial^{2}}{\partial s^{2}}G(B^{w}_{r},1)\bp_{\delta_{B^{w}_{r}}}\left[\langle f,Z_{t-r}\rangle\right]^{2}dr\right]\nonumber\\
&=&\frac{1}{w(x)}P^{\beta}_{t}(w f^{2})(x)+\frac{1}{w(x)}\int_{0}^{t}P^{\beta}_{s}\left[2\gamma\left(P^{\beta}_{t-s}(w f)\right)^{2}\right](x)ds.\label{varianceZ}
\end{eqnarray}
One can easily show by \eqref{meanZ} that
\begin{equation}\label{martingale-X-Z}
W^{h/w}_{t}(Z):=\e^{-\lambda_{1}t}\langle \frac{h}{w},Z_{t}\rangle,
\quad t\ge 0,
\end{equation}
is a nonnegative $\bp_{\delta_{x}}$-martingale for every $x\in\R$, and a nonnegative $\p_{\mu}$-martingale for every $\mu\in\mf$.
We use $W^{h/w}_{\infty}(Z)$ to denote the martingale limit. It is proved by \cite[Proposition 1.1]{EKW} that
\begin{equation}\label{martingale-limit-X-Z}
W^{h/w}_{\infty}(Z)=W^{h}_{\infty}(X)\quad\p_{\mu}\mbox{-a.s.}
\end{equation}
for all $\mu\in\mf$.

\begin{lemma}\label{lem4.4}
\begin{description}
\item{(i)} There exist $T_{1}>0$ and $\theta_{i}(t)$ ($i=1,2$) such that for $t\ge  T_{1}$,
\begin{equation}\label{lem4.4.1}
\theta_{1}(t)\le \frac{\bp_{\delta_{x(t)}}\left[Z_{t-s(t)}(A+R(t))\right]}{w_{-}\eta_{-}(A)\frac{h(x(t))}{w(x(t))}\e^{\lambda_{1}(t-s(t))-\sqrt{2\lambda_{1}}R(t)}}\le \theta_{2}(t),
\end{equation}
where
 $\eta_{-}(dx)=C_{-}\e^{-\sqrt{2\lambda_{1}}x}dx$, and
for $i=1,2$, $\theta_{i}(t)\to 1$ as $t\to+\infty$.

\item{(ii)} There exist $T_{2},C>0$ such that for $t\ge T_{2}$,
\begin{multline}\label{lem4.4.2}
\bp_{\delta_{x(t)}}\left[Z_{t-s(t)}(A+R(t))\right]\le \bp_{\delta_{x(t)}}\left[Z_{t-s(t)}(A+R(t))^{2}\right]\\
\le \bp_{\delta_{x(t)}}\left[Z_{t-s(t)}(A+R(t))\right]+C\frac{h(x(t))}{w(x(t))}\e^{2\lambda_{1}(t-s(t))-2\sqrt{2\lambda_{1}}R(t)}.
\end{multline}
\end{description}
\end{lemma}

\proof (i) We have
$$\bp_{\delta_{x(t)}}\left[Z_{t-s(t)}(A+R(t))\right]=\frac{1}{w(x(t))}P^{\beta}_{t-s(t)}\left(w\1_{A+R(t)}\right)(x(t)).$$
Note that for $t$ large enough such
that $R(t)+(\inf A\wedge 0)\ge M$, $w(x)=w_{-}$ for all $x\in A+R(t)$. It follows that
$$\bp_{\delta_{x(t)}}\left[Z_{t-s(t)}(A+R(t))\right]=\frac{w_{-}}{w(x(t))}P^{\beta}_{t-s(t)}\1_{A+R(t)}(x(t)).$$
Thus \eqref{lem4.4.1} follows immediately from Lemma \ref{lem:esti2}.

\medskip

(ii) The first inequality of \eqref{lem4.4.2} is obvious
since, by \eqref{varianceZ},
\begin{eqnarray}\label{lem4.4.3}
\bp_{\delta_{x(t)}}\left[Z_{t-s(t)}(A+R(t))^{2}\right]&=&\bp_{\delta_{x(t)}}\left[Z_{t-s(t)}(A+R(t))\right]\nonumber\\
&&+\frac{2}{w(x(t))}\int_{0}^{t-s(t)}P^{\beta}_{r}\left[\gamma P^{\beta}_{t-s(t)-r}\left(w\1_{A+R(t)}\right)^{2}\right](x(t))dr.
\end{eqnarray}
Suppose supp$\gamma\subset [-k,k]$ for some $0<k<+\infty$.
Let $\sigma_{k}$ be the first hitting time of $[-k,k]$ by the Brownian motion.
Noting that $s(t)\le a t<t-1$ for $t$ sufficiently large,
we have
\begin{eqnarray}
&&\int_{0}^{t-s(t)}P^{\beta}_{r}\left[\gamma \left(P^{\beta}_{t-s(t)-r}\1_{A+R(t)}\right)^{2}\right](x(t))dr\nonumber\\
&=&\Bp_{x(t)}\left[\int_{0}^{t-s(t)}e_{\beta}(r)\gamma(B_{r})\left(P^{\beta}_{t-s(t)-r}\1_{A+R(t)}(B_{r})\right)^{2}dr\right]\nonumber\\
&=&\Bp_{x(t)}\left[\int_{\sigma_{k}}^{t-s(t)}e_{\beta}(r)\gamma(B_{r})\left(P^{\beta}_{t-s(t)-r}\1_{A+R(t)}(B_{r})\right)^{2}dr;\sigma_{k}\le t-s(t)\right]\nonumber\\
&=&\Bp_{x(t)}\left[e_{\beta}(\sigma_{k})\left.\Bp_{B_{u}}\left[\int_{0}^{t-s(t)-u}e_{\beta}(r)\gamma(B_{r})\left(P^{\beta}_{t-s(t)-u-r}\1_{A+R(t)}(B_{r})\right)^{2}dr\right]\right|_{u=\sigma_{k}};\sigma_{k}\le t-s(t)\right]\nonumber\\
&\le&c_{1}\Bp_{x(t)}\left[e_{\beta}(\sigma_{k})\frac{h(B_{\sigma_{k}})}{\inf_{x\in [-k,k]}h(x)}\e^{2\lambda_{1}(t-s(t)-\sigma_{k})-2\sqrt{2\lambda_{1}}R(t)};\sigma_{k}\le t-s(t)\right]\nonumber\\
&=&c_{2}\e^{2\lambda_{1}(t-s(t))-2\sqrt{2\lambda_{1}}R(t)}\Bp_{x(t)}\left[e_{\beta}(\sigma_{k})h(B_{\sigma_{k}})\e^{-2\lambda_{1}\sigma_{k}};\sigma_{k}\le t-s(t)\right].\label{lem4.4.4}
\end{eqnarray}
The above inequality
follows from Lemma \ref{lem:estimate}(ii).
Since $\e^{-\lambda_{1}t}\e_{\beta}(t)h(B_{t})$ is a martingale, by the optional stopping theorem, the last term in \eqref{lem4.4.4} is no larger than $h(x(t))$. So we get that
$$\int_{0}^{t-s(t)}P^{\beta}_{r}\left[\gamma \left(P^{\beta}_{t-s(t)-r}\1_{A+R(t)}\right)^{2}\right](x(t))dr\le c_{3}\e^{2\lambda_{1}(t-s(t))-2\sqrt{2\lambda_{1}}R(t)}h(x(t)).$$
We also note that for $t$ large enough, $w(x)=w_{-}$ for all $x\in A+R(t)$. Thus
\begin{eqnarray*}
\int_{0}^{t-s(t)}P^{\beta}_{r}\left[\gamma P^{\beta}_{t-s(t)-r}\left(w\1_{A+R(t)}\right)^{2}\right](x(t))dr&=&w^{2}_{-}\int_{0}^{t-s(t)}P^{\beta}_{r}\left[\gamma \left(P^{\beta}_{t-s(t)-r}\1_{A+R(t)}\right)^{2}\right](x(t))dr\nonumber\\
&\le&c_{3}w_{-}^{2}\e^{2\lambda_{1}(t-s(t))-2\sqrt{2\lambda_{1}}R(t)}h(x(t)).
\end{eqnarray*}
Putting this back to \eqref{lem4.4.3} we get \eqref{lem4.4.2}.\qed

\bigskip

\begin{lemma}\label{lem4.5}
Assume that $\delta=\sqrt{\lambda_{1}/2}$ and that
 $\lambda_{1}s(t)+\sqrt{2\lambda_{1}}a(t)\to +\infty$ as $t\to+\infty$. Then there exist $C,T>0$ and $\theta_{i}(t)$ ($i=4,5,6,7$) such that
for $t\ge T$,
\begin{eqnarray}
 \bp_{\delta_{x(t)}}\left(Z_{t-s(t)}(A+R(t))=0\right)&\le& 1-\theta_{4}(t)w_{-}\eta_{-}(A)\frac{h(x(t))}{w(x(t))}\Theta(t),\label{lem4.5.1}\\
 \bp_{\delta_{x(t)}}\left(Z_{t-s(t)}(A+R(t))=0\right)&\ge& 1-\theta_{5}(t)w_{-}\eta_{-}(A)\frac{h(x(t))}{w(x(t))}\Theta(t),\label{lem4.5.2}\\
 \bp_{\delta_{x(t)}}\left(Z_{t-s(t)}(A+R(t))=1\right)&\le& \theta_{6}(t)w_{-}\eta_{-}(A)\frac{h(x(t))}{w(x(t))}\Theta(t),\label{lem4.5.3}\\
 \bp_{\delta_{x(t)}}\left(Z_{t-s(t)}(A+R(t))=1\right)&\ge& \theta_{7}(t)w_{-}\eta_{-}(A)\frac{h(x(t))}{w(x(t))}\Theta(t),\label{lem4.5.4}\\
 \bp_{\delta_{x(t)}}\left(Z_{t-s(t)}(A+R(t))\ge 2\right)&\le&C\frac{h(x(t))}{w(x(t))}\Theta^{2}(t).\label{lem4.5.5}
 \end{eqnarray}
where $\Theta(t)=\e^{-\lambda_{1}s(t)-\sqrt{2\lambda_{1}}a(t)}$ and $\theta_{i}(t)\to 1$ as $t\to+\infty$ for $i=4,5,6,7$ .
\end{lemma}

\proof We note that if $Z$ is an integer-valued random variable then
\begin{equation}\label{lem4.5.6}
\frac{E[Z]^{2}}{E[Z^{2}]}\le P(Z>0)=P(Z\ge 1)\le E[Z],
\end{equation}
and
\begin{equation}\label{lem4.5.7}
P(Z\ge 2)\le E[Z(Z-1);Z\ge 2]=E[Z(Z-1)]=E[Z^{2}]-E[Z].
\end{equation}
It is easy to see that \eqref{lem4.5.5} follows immediately from \eqref{lem4.5.7} and Lemma \ref{lem4.4}(ii).
Since $1-\bp_{\delta_{x(t)}}\left(Z_{t-s(t)}(A+R(t))=0\right)=\bp_{\delta_{x(t)}}\left(Z_{t-s(t)}(A+R(t))>0\right)$, we have by \eqref{lem4.5.6} and Lemma \ref{lem4.4} that for $t$ large enough,
$$1-\bp_{\delta_{x(t)}}\left(Z_{t-s(t)}(A+R(t))=0\right)\le \bp_{\delta_{x(t)}}\left[Z_{t-s(t)}(A+R(t))\right] \le \theta_{2}(t)w_{-}\eta_{-}(A)\frac{h(x(t))}{w(x(t))}\Theta(t),$$
and
\begin{eqnarray*}
1-\bp_{\delta_{x(t)}}\left(Z_{t-s(t)}(A+R(t))=0\right)&\ge&\frac{\bp_{\delta_{x(t)}}\left[Z_{t-s(t)}(A+R(t))\right]^{2}}{\bp_{\delta_{x(t)}}\left[Z_{t-s(t)}(A+R(t))^{2}\right]}\\
&\ge&\frac{\bp_{\delta_{x(t)}}\left[Z_{t-s(t)}(A+R(t))\right]^{2}}{\bp_{\delta_{x(t)}}\left[Z_{t-s(t)}(A+R(t))\right]+C\frac{h(x(t))}{w(x(t))}\Theta^{2}(t)}\\
&\ge&\frac{\left[\theta_{1}(t)w_{-}\eta_{-}(A)\frac{h(x(t))}{w(x(t))}\Theta(t)\right]^{2}}{\theta_{2}(t)w_{-}\eta_{-}(A)\frac{h(x(t))}{w(x(t))}\Theta(t)+C\frac{h(x(t))}{w(x(t))}\Theta^{2}(t)}\\
&=&\frac{\theta_{1}(t)^{2}}{\theta_{2}(t)+Cw^{-1}_{-}\eta^{-1}_{-}(A)\Theta(t)}\,w_{-}\eta_{-}(A)\frac{h(x(t))}{w(x(t))}\Theta(t).
\end{eqnarray*}
Since $\theta_{i}(t)\to 1$ for $i=1,2$ and $\Theta(t)\to 0$ as $t\to+\infty$, $\frac{\theta_{1}(t)^{2}}{\theta_{2}(t)+Cw^{-1}_{-}\eta^{-1}_{-}(A)\Theta(t)}\to 1$ as $t\to+\infty$. Hence we prove \eqref{lem4.5.1} and \eqref{lem4.5.2}.

We note that
$$\bp_{\delta_{x(t)}}\left(Z_{t-s(t)}(A+R(t))=1\right)=1-\bp_{\delta_{x(t)}}\left(Z_{t-s(t)}(A+R(t))=0\right)-\bp_{\delta_{x(t)}}\left(Z_{t-s(t)}(A+R(t))\ge 2\right).$$
Thus \eqref{lem4.5.3} and \eqref{lem4.5.4} follow immediately from \eqref{lem4.5.1}, \eqref{lem4.5.2} and \eqref{lem4.5.5}.\qed

\bigskip

\begin{lemma}\label{lem4.6}
\begin{description}
\item{(i)} For every $\delta\in (\sqrt{\lambda_{1}/2},\sqrt{2\lambda_{1}})$ and $x\in\R$,
$$\lim_{t\to+\infty}\bp_{\delta_{x}}\left(\max_{u\in Z_{t}}|z_{u}(t)|<\delta t\right)=1.$$
\item{(ii)} For every $x\in \R$,
$$\liminf_{t\to+\infty}\e^{-\lambda_{1}t}\|Z_{t}\|>0\quad\bp_{\delta_{x}}\mbox{-a.s.}$$
\end{description}
\end{lemma}

\proof (i) We have
\begin{eqnarray}
\bp_{\delta_{x}}\left(\max_{u\in Z_{t}}|z_{u}(t)|<\delta t\right)&=&\bp_{\delta_{x}}\left(Z_{t}((-\delta t,\delta t)^{c})=0\right)\nonumber\\
&=&1-\bp_{\delta_{x}}\left(Z_{t}((-\delta t,\delta t)^{c})\ge 1\right)\nonumber\\
&\ge&1-\bp_{\delta_{x}}\left[Z_{t}((-\delta t,\delta t)^{c})\right].\nonumber
\end{eqnarray}
So it suffices to show that
\begin{equation}\label{lem4.6.1}
\lim_{t\to+\infty}\bp_{\delta_{x}}\left[Z_{t}((-\delta t,\delta t)^{c})\right]=0.
\end{equation}
Note that for $t$ large enough such that $\delta t>M$,
\begin{eqnarray}
\bp_{\delta_{x}}\left[Z_{t}((-\delta t,\delta t)^{c})\right]&=&\frac{1}{w(x)}P^{\beta}_{t}\left(w\1_{(-\delta t,\delta t)^{c}}\right)(x)\nonumber\\
&\le&\frac{w_{+}\vee w_{-}}{w(x)}P^{\beta}_{t}\1_{(-\delta t,\delta t)^{c}}(x)=
 \frac{w_{+}\vee w_{-}}{w(x)}
 \pi^{\delta t}_{t}(x).\nonumber
\end{eqnarray}
It follows by Lemma \ref{lem1} that $\lim_{t\to+\infty}\log\pi^{\delta t}_{t}(x)/t=-\Lambda_{\delta}<0$ for $\delta\in (\sqrt{\lambda_{1}/2},\sqrt{2\lambda_{1}})$.
So we get $\pi^{\delta t}_{t}(x)\to 0$
as $t\to+\infty$ and \eqref{lem4.6.1} follows immediately.

\medskip

(ii) We note that for $t\ge 0$,
$$\e^{-\lambda_{1}t}\|Z_{t}\|=\e^{-\lambda_{1}t}\Big\langle \frac{w}{h}\cdot\frac{h}{w},Z_{t}\Big\rangle\ge \|\frac{h}{w}\|^{-1}_{\infty}W^{h/w}_{t}(Z),$$
where $\|h/w\|_{\infty}<+\infty$. So it suffices to show that for every $x\in\R$,
\begin{equation}\label{lem4.6.2}
\bp_{\delta_{x}}\left(W^{h/w}_{\infty}(Z)>0\right)=1.
\end{equation}
We have
\begin{equation}\label{lem4.6.3}
\p_{\delta_{x}}\left(W^{h/w}_{\infty}(Z)=0\right)=\p_{\delta_{x}}\left[\prod_{u\in Z_{0}}\bp_{\delta_{z_{u}(0)}}\left(W^{h/w}_{\infty}(Z)=0\right)\right]=\e^{-w(x)\left(1-\bp_{\delta_{x}}\left(W^{h/w}_{\infty}(Z)=0\right)\right)}.
\end{equation}
The final equality is because $(Z_{0},\p_{\delta_{x}})$ is a Poisson point process with intensity $w\delta_{x}$.
On the other hand, by \eqref{martingale-limit-X-Z} we have
$$\p_{\delta_{x}}\left(W^{h/w}_{\infty}(Z)=0\right)=\p_{\delta_{x}}\left(W^{h}_{\infty}(X)=0\right)=\e^{-w(x)}.$$
Combining this with \eqref{lem4.6.3}, we get that $\bp_{\delta_{x}}\left(W^{h/w}_{\infty}(Z)=0\right)=0$ and \eqref{lem4.6.2} follows immediately.\qed

\bigskip

\begin{proposition}\label{prop4}
For every $x\in\R$, $((Z_{t}\pm\sqrt{\lambda_{1}/2}\,t)_{t\ge 0},\bp_{\delta_{x}})$ converges in distribution to a Cox process directed by
$w_{\pm}W^{h/w}_{\infty}(Z)\eta_{\pm}(dx)$,
where $W^{h/w}_{\infty}(Z)$ is the martingale limit of
$((W^{h/w}_{t}(Z))_{t\geq 0}, \bp_{\delta_{x}})$.
\end{proposition}

\proof
Take $R(t)=\sqrt{\lambda_{1}/2}\,t$ and fix a function $s(\cdot)$ such that $s(t)\to +\infty$ and $s(t)=o(t)$ as $t\to+\infty$.
For notational simplicity, in  the proof we shall write $s(t)$ as $s$.
We only consider $Z_{t}-\sqrt{\lambda_{1}/2}\,t$. The result for $Z_{t}+\sqrt{\lambda_{1}/2}\,t$ can be proved similarly.

It follows from Lemma \ref{lem4.4}(i)
that for any Borel set $A$ with $\inf A>-\infty$,
\begin{equation}\label{cond:ii}
\bp_{\delta_{x}}\left[Z_{t}(A+R(t))\right]\sim w_{-}\eta_{-}(A)\frac{h(x)}{w(x)}\mbox{ as } t\to+\infty.
\end{equation}
We
only need to show that condition (i) of Lemma \ref{lemA1} is satisfied by any subsequence of $\{Z_{t}-\sqrt{\lambda_{1}/2}\,t:t\ge 0\}$.

Take $m\in\mathbb{N}$, $k_{1},\cdots,k_{m}\in\mathbb{Z}^{+}$ and mutually disjoint Borel sets $A_{1},\cdots,A_{m}$ in $\R$
with $\inf A_{i}>-\infty$ for $i=1,\cdots, m$.
Put
$k:=k_{1}+\cdots+k_{m}$ and $A:=\bigcup_{i=1}^{m}A_{i}$. Let $\mathcal{G}_{s}$ be the $\sigma$-field generated by $\{Z_{r}:r\in [0,s]\}$.
It suffices to show that
\begin{equation}\label{prop4.5}
\bp_{\delta_{x}}\left(\bigcap_{i=1}^{m}\left\{Z_{t}(A_{i}+R(t))=k_{i}\right\}\,|\,\mathcal{G}_{s}\right)\to \e^{-w_{-}W^{h/w}_{\infty}(Z)\sum_{i=1}^{m}\eta_{-}(A_{i})}\prod_{i=1}^{m}\frac{\left(w_{-}W^{h/w}_{\infty}(Z)\eta_{-}(A_{i})\right)^{k_{i}}}{k_{i}!}.
\end{equation}
in probability as $t\to+\infty$.
\footnote{Actually \eqref{prop4.5} is a bit stronger than what one needs for the proof of Proposition \ref{prop4}. The proof can be shortened by applying \cite[Proposition 16.17]{Kallenberg2}. In fact by the aforementioned result,
one only needs to show that (i) \eqref{cond:ii}
holds for all relatively compact sets $A\subset \R$, and (ii) $\lim_{t\to +\infty}\bp_{\delta_{x}}\left(Z_{t}(A+R(t))=0\right)=\bp_{\delta_{x}}\left[\exp\{-w_{-}W^{h/w}_{\infty}(Z)\eta_{-}(A)\}\right]$
for all compact sets $A$. However, since \eqref{prop4.5} further yields the limit of the order statistics of $Z_{t}$ (see Proposition \ref{prop6} and the remark below), we present it here for the sake of being more self-contained.}
For $u\in Z_{s}$, let $Z^{(u)}_{t-s}$ be the point process of the locations of the particles alive at time $t$ whose ancestor is $u$.
Take a constant $\kappa >\sqrt{\lambda_{1}/2}$. Define
$$\mathcal{E}^{1}_{t}:=\left\{\max_{u\in Z_{s}}|z_{u}(s)|\le \kappa s,\ \|Z_{s}\|\ge k\right\}\mbox{ and }\mathcal{E}^{2}_{t}:=\left\{Z^{(u)}_{t-s}(A+R(t))\le 1\quad\forall u\in Z_{s}\right\}.$$
It follows from Lemma \ref{lem4.6} that
$\bp_{\delta_{x}}\left(\mathcal{E}^{1}_{t}\right)\to 1$ as $t\to+\infty$.
Since $\mathcal{E}^{1}_{t}\in \mathcal{G}_{s}$, we get
\begin{equation}\label{prop4.2}
\bp_{\delta_{x}}\left((\mathcal{E}^{1}_{t})^{c}\,|\,\mathcal{G}_{s}\right)=\1_{(\mathcal{E}^{1}_{t})^{c}}\to 0
\mbox{ as }t\to+\infty\quad\mbox{ in probability.}
\end{equation}
On the event $\mathcal{E}^{1}_{t}$ we have
\begin{eqnarray}
\bp_{\delta_{x}}\left((\mathcal{E}^{2}_{t})^{c}\,|\,\mathcal{G}_{s}\right)
&=&\bp_{\delta_{x}}\left(Z^{(u)}_{t-s}(A+R(t))\ge 2\mbox{ for some } u\in Z_{s}\,|\,\mathcal{G}_{s}\right)\nonumber\\
&\le&\sum_{u\in Z_{s}}\bp_{\delta_{z_{u}(s)}}\left(Z^{(u)}_{t-s}(A+R(t))\ge 2\right).\label{prop4.1}
\end{eqnarray}
By Lemma \ref{lem4.5}, for $t$ large enough, on the event $\mathcal{E}^{1}_{t}$, $\bp_{\delta_{z_{u}(s)}}\left(Z^{(u)}_{t-s}(A+R(t))\ge 2\right)\le c_{1}\frac{h(z_{u}(s))}{w(z_{u}(s))}\e^{-2\lambda_{1}s}$. Hence we get by \eqref{prop4.1} that on $\mathcal{E}^{1}_{t}$
\begin{equation}\nonumber
\bp_{\delta_{x}}\left((\mathcal{E}^{2}_{t})^{c}\,|\,\mathcal{G}_{s}\right)
\le c_{1}\sum_{u\in Z_{s}}\frac{h(z_{u}(s))}{w(z_{u}(s))}\e^{-2\lambda_{1}s}
\le c_{1}\e^{-\lambda_{1}s}W^{h/w}_{s}(Z).
\end{equation}
This yields $\1_{\mathcal{E}^{1}_{t}}\bp_{\delta_{x}}\left((\mathcal{E}^{2}_{t})^{c}\,|\,\mathcal{G}_{s}\right)\to 0$ $\bp_{\delta_{x}}$-a.s. Consequently by \eqref{prop4.2} we have
\begin{equation}\label{prop4.3}
\bp_{\delta_{x}}\left((\mathcal{E}^{2}_{t})^{c}\,|\,\mathcal{G}_{s}\right)\to 0\quad \mbox{ in probability as }t\to+\infty.
\end{equation}
By \eqref{prop4.2} and \eqref{prop4.3}, we have
\begin{equation}
\bp_{\delta_{x}}\left(\bigcap_{i=1}^{m}\left\{Z_{t}(A_{i}+R(t))=k_{i}\right\}\,|\,\mathcal{G}_{s}\right)
=\bp_{\delta_{x}}\left(\bigcap_{i=1}^{m}\left\{Z_{t}(A_{i}+R(t))=k_{i}\right\},\mathcal{E}^{1}_{t},\mathcal{E}^{2}_{t}\,|\,\mathcal{G}_{s}\right)+\epsilon^{1}_{t}
\end{equation}
for some $\epsilon^{1}_{t}\to 0$ in probability.

We note that on the event $\mathcal{E}^{2}_{t}$, $\{Z^{(u)}_{t-s}(A+R(t)):u\in Z_{s}\}$ are Bernoulli random variables. So we have
\begin{eqnarray}\label{prop4.4}
&&\bp_{\delta_{x}}\left(\bigcap_{i=1}^{m}\left\{Z_{t}(A_{i}+R(t))=k_{i}\right\},\mathcal{E}^{1}_{t},\mathcal{E}^{2}_{t}\,|\,\mathcal{G}_{s}\right)\nonumber\\
&=&\frac{1}{k_{1}!\cdots k_{m}!}\bp_{\delta_{x}}\Big(\bigcup_{(u_{1},\cdots,u_{k})\subset Z_{s}}\big\{Z^{(u_{1})}_{t-s}(A_{1}+R(t))=1,\cdots,Z^{(u_{k})}_{t-s}(A_{m}+R(t))=1,\nonumber\\ &&\quad\quad\quad\quad\quad Z^{(u)}_{t-s}(A+R(t))=0\quad\forall u\in Z_{s}\setminus\{u_{1},\cdots,u_{k}\}\big\},\mathcal{E}^{1}_{t},\mathcal{E}^{2}_{t}\,|\,\mathcal{G}_{s}\Big),
\end{eqnarray}
where $(u_{1},\cdots,u_{k})\subset Z_{s}$ is the union over all $k$-permutations of $Z_{s}$. By \eqref{prop4.3} and the fact that $\mathcal{E}^{1}_{t}\in\mathcal{G}_{s}$, the conditional probability in the right hand side of \eqref{prop4.4} equals
\begin{eqnarray}
&&\1_{\mathcal{E}^{1}_{t}}\bp_{\delta_{x}}\Big(\bigcup_{(u_{1},\cdots,u_{k})\subset Z_{s}}\big\{Z^{(u_{1})}_{t-s}(A_{1}+R(t))=1,\cdots,Z^{(u_{k})}_{t-s}(A_{m}+R(t))=1,\nonumber\\ &&\quad\quad\quad\quad\quad Z^{(u)}_{t-s}(A+R(t))=0\quad\forall u\in Z_{s}\setminus\{u_{1},\cdots,u_{k}\}\big\}\,|\,\mathcal{G}_{s}\Big)+\epsilon^{2}_{t},
\end{eqnarray}
where $\epsilon^{2}_{t}\to 0$ in probability.
Since $\bigcup_{(u_{1},\cdots,u_{k})\subset Z_{s}}\{\cdots\}$ is
an union of mutually-disjoint events, we have
\begin{eqnarray}
&&\bp_{\delta_{x}}\Big(\bigcup_{(u_{1},\cdots,u_{k})\subset Z_{s}}\big\{Z^{(u_{1})}_{t-s}(A_{1}+R(t))=1,\cdots,Z^{(u_{k})}_{t-s}(A_{m}+R(t))=1,\nonumber\\
&&\quad\quad\quad\quad\quad Z^{(u)}_{t-s}(A+R(t))=0\quad\forall u\in Z_{s}\setminus\{u_{1},\cdots,u_{k}\}\big\}\,|\,\mathcal{G}_{s}\Big)\nonumber\\
&=&\sum_{(u_{1},\cdots,u_{k})\subset Z_{s}}\bp_{\delta_{x}}\big(Z^{(u_{1})}_{t-s}(A_{1}+R(t))=1,\cdots,Z^{(u_{k})}_{t-s}(A_{m}+R(t))=1,\nonumber\\
&&\quad\quad\quad\quad\quad Z^{(u)}_{t-s}(A+R(t))=0\quad\forall u\in Z_{s}\setminus\{u_{1},\cdots,u_{k}\}\,|\,\mathcal{G}_{s}\big)\nonumber\\
&=&\sum_{(u_{1},\cdots,u_{k})\subset Z_{s}}\bp_{\delta_{z_{u_{1}}(s)}}\big(Z_{t-s}(A_{1}+R(t))=1\big)\times \cdots \times
\bp_{\delta_{z_{u_{k}}(s)}}\big(Z_{t-s}(A_{m}+R(t))=1\big)\nonumber\\
&&\quad\quad\quad\quad \times \prod_{u\in Z_{s}\setminus\{u_{1},\cdots,u_{k}\}}\bp_{\delta_{z_{u}(s)}}\big( Z_{t-s}(A+R(t))=0\big)\nonumber\\
&=&\Big[\prod_{u\in Z_{s}}\bp_{\delta_{z_{u}(s)}}\big( Z_{t-s}(A+R(t))=0\big)\Big]\times\nonumber\\
&&\Big[\sum_{(u_{1},\cdots,u_{k})\subset Z_{s}}\frac{\bp_{\delta_{z_{u_{1}}(s)}}\big(Z_{t-s}(A_{1}+R(t))=1\big)}{\bp_{\delta_{z_{u_{1}}(s)}}\big(Z_{t-s}(A_{1}+R(t))=0\big)}\times \cdots \times
\frac{\bp_{\delta_{z_{u_{k}}(s)}}\big(Z_{t-s}(A_{m}+R(t))=1\big)}{\bp_{\delta_{z_{u_{k}}(s)}}\big(Z_{t-s}(A_{1}+R(t))=0\big)}\Big].\nonumber
\end{eqnarray}

The second equality follows the Markov branching property. So far we have proved that
\begin{eqnarray}
&&\bp_{\delta_{x}}\left(\bigcap_{i=1}^{m}\left\{Z_{t}(A_{i}+R(t))=k_{i}\right\}\,|\,\mathcal{G}_{s}\right)\nonumber\\
&=&\epsilon^{1}_{t}+\epsilon^{2}_{t}+\frac{1}{k_{1}!\cdots k_{m}!}\1_{\mathcal{E}^{1}_{t}}\Big[\prod_{u\in Z_{s}}\bp_{\delta_{z_{u}(s)}}\big( Z_{t-s}(A+R(t))=0\big)\Big]\times\nonumber\\
&&\Big[\sum_{(u_{1},\cdots,u_{k})\subset Z_{s}}\frac{\bp_{\delta_{z_{u_{1}}(s)}}\big(Z_{t-s}(A_{1}+R(t))=1\big)}{\bp_{\delta_{z_{u_{1}}(s)}}\big(Z_{t-s}(A_{1}+R(t))=0\big)}\times \cdots \times
\frac{\bp_{\delta_{z_{u_{k}}(s)}}\big(Z_{t-s}(A_{m}+R(t))=1\big)}{\bp_{\delta_{z_{u_{k}}(s)}}\big(Z_{t-s}(A_{1}+R(t))=0\big)}\Big].\nonumber
\end{eqnarray}
Hence to prove \eqref{prop4.5}, it suffices to prove that
\begin{equation}\label{prop4.6}
\lim_{t\to+\infty}\1_{\mathcal{E}^{1}_{t}}\prod_{u\in Z_{s}}\bp_{\delta_{z_{u}(s)}}\big( Z_{t-s}(A+R(t))=0\big)=\e^{-w_{-}W^{h/w}_{\infty}(Z)\eta_{-}(A)}\quad\mbox{ in probability},
\end{equation}
and
\begin{eqnarray}\label{prop4.7}
&&\lim_{t\to+\infty}\1_{\mathcal{E}^{1}_{t}}\sum_{(u_{1},\cdots,u_{k})\subset Z_{s}}\frac{\bp_{\delta_{z_{u_{1}}(s)}}\big(Z_{t-s}(A_{1}+R(t))=1\big)}{\bp_{\delta_{z_{u_{1}}(s)}}\big(Z_{t-s}(A_{1}+R(t))=0\big)}\times \cdots \times
\frac{\bp_{\delta_{z_{u_{k}}(s)}}\big(Z_{t-s}(A_{m}+R(t))=1\big)}{\bp_{\delta_{z_{u_{k}}(s)}}\big(Z_{t-s}(A_{1}+R(t))=0\big)}\nonumber\\
&=&\prod_{i=1}^{m}\left(w_{-}W^{h/w}_{\infty}(Z)\eta_{-}(A_{i})\right)^{k_{i}}\quad\mbox{ in probability.}
\end{eqnarray}

(i) We first prove \eqref{prop4.6}.
It follows from \eqref{lem4.5.1} that for $t$ large enough, on the event $\mathcal{E}^{1}_{t}$,
\begin{eqnarray}\label{prop4.8}
\prod_{u\in Z_{s}}\bp_{\delta_{z_{u}(s)}}\big( Z_{t-s}(A+R(t))=0\big)
&\le&\prod_{u\in Z_{s}}\left(1-\theta_{4}(t)w_{-}\eta_{-}(A)\frac{h(z_{u}(s))}{w(z_{u}(s))}\e^{-\lambda_{1}s}\right)\nonumber\\
&\le&\prod_{u\in Z_{s}}\exp\{-\theta_{4}(t)w_{-}\eta_{-}(A)\frac{h(z_{u}(s))}{w(z_{u}(s))}\e^{-\lambda_{1}s}\}\nonumber\\
&=&\exp\{-\theta_{4}(t)w_{-}\eta_{-}(A)W^{h/w}_{s}(Z)\}.
\end{eqnarray}
The second inequality is from the fact that $1-x\le \e^{-x}$ for all $x\ge 0$.
For the lower bound, it follows from  \eqref{lem4.5.2} that for $t$ large enough, on $\mathcal{E}^{1}_{t}$
\begin{eqnarray}
&&\prod_{u\in Z_{s}}\bp_{\delta_{z_{u}(s)}}\big( Z_{t-s}(A+R(t))=0\big)\nonumber\\
&\ge&\prod_{u\in Z_{s}}\left(1-\theta_{5}(t)w_{-}\eta_{-}(A)\frac{h(z_{u}(s))}{w(z_{u}(s))}\e^{-\lambda_{1}s}\right)\nonumber\\
&=&\exp\left\{\sum_{u\in Z_{s}}\log\left(1-\theta_{5}(t)w_{-}\eta_{-}(A)\frac{h(z_{u}(s))}{w(z_{u}(s))}\e^{-\lambda_{1}s}\right)\right\}.\label{prop4.9}
\end{eqnarray}
Note that $c:=\sup_{y\in\R}h(y)/w(y)\le \|h\|_{\infty}/\inf_{y\in\R}w(y)<+\infty$. Using the fact that
$$\log(1-x)\ge \frac{\log (1-x^{*})}{x^{*}}\,x\quad\forall x^{*}\in (0,1),\ x\in [0,x^{*}],$$
we get by \eqref{prop4.9} that on $\mathcal{E}^{1}_{t}$
\begin{eqnarray}
&&\prod_{u\in Z_{s}}\bp_{\delta_{z_{u}(s)}}\big( Z_{t-s}(A+R(t))=0\big)\nonumber\\
&\ge&\exp\Big\{\frac{\log\left(1-\theta_{5}(t)w_{-}\eta_{-}(A)c\e^{-\lambda_{1}s}\right)}{\theta_{5}(t)w_{-}\eta_{-}(A)c\e^{-\lambda_{1}s}}\,\sum_{u\in Z_{s}}\theta_{5}(t)w_{-}\eta_{-}(A)\frac{h(z_{u}(s))}{w(z_{u}(s))}\e^{-\lambda_{1}s}\Big\}\nonumber\\
&=&\exp\Big\{\frac{\log\left(1-\theta_{5}(t)w_{-}\eta_{-}(A)c\e^{-\lambda_{1}s}\right)}{\theta_{5}(t)w_{-}\eta_{-}(A)c\e^{-\lambda_{1}s}}\,
\theta_{5}(t)w_{-}\eta_{-}(A)W^{h/w}_{s}(Z)\Big\}.\label{prop4.10}
\end{eqnarray}
It is easy to see that the final terms of \eqref{prop4.8} and \eqref{prop4.10}
converges to $\exp\left\{-w_{-}\eta_{-}(A)W^{h/w}_{\infty}(Z)\right\}$ almost surely.
Thus \eqref{prop4.6} follows immediately.

\medskip

(ii) Now we prove \eqref{prop4.7}.
We use $\theta^{(j)}_{i}(t)$ to denote the functions $\theta_{i}(t)$ in Lemma \ref{lem4.5} corresponding to the set $A_{j}$. It follows by \eqref{lem4.5.2} and \eqref{lem4.5.3} that, on $\mathcal{E}^{1}_{t}$,
\begin{eqnarray}\label{prop4.11}
&&\sum_{(u_{1},\cdots,u_{k})\subset Z_{s}}\frac{\bp_{\delta_{z_{u_{1}}(s)}}\big(Z_{t-s}(A_{1}+R(t))=1\big)}{\bp_{\delta_{z_{u_{1}}(s)}}\big(Z_{t-s}(A_{1}+R(t))=0\big)}\times \cdots \times
\frac{\bp_{\delta_{z_{u_{k}}(s)}}\big(Z_{t-s}(A_{m}+R(t))=1\big)}{\bp_{\delta_{z_{u_{k}}(s)}}\big(Z_{t-s}(A_{1}+R(t))=0\big)}\nonumber\\
&\le&\sum_{(u_{1},\cdots,u_{k})\subset Z_{s}}\frac{\theta^{(1)}_{6}(t)w_{-}\eta_{-}(A_{1})\frac{h(z_{u_{1}}(s))}{w(z_{u_{1}}(s))}\e^{-\lambda_{1}s}}{1-\theta_{5}(t)w_{-}\eta_{-}(A)\frac{h(z_{u_{1}}(s))}{w(z_{u_{1}}(s))}\e^{-\lambda_{1}s}}
\times\cdots\times \frac{\theta^{(m)}_{6}(t)w_{-}\eta_{-}(A_{m})\frac{h(z_{u_{k}}(s))}{w(z_{u_{k}}(s))}\e^{-\lambda_{1}s}}{1-\theta_{5}(t)w_{-}\eta_{-}(A)\frac{h(z_{u_{k}}(s))}{w(z_{u_{k}}(s))}\e^{-\lambda_{1}s}}\nonumber\\
&\le&\frac{\prod_{i=1}^{m}\theta^{(i)}_{6}(t)^{k_{i}}}{\left(1-\theta_{5}(t)w_{-}\eta_{-}(A)c\e^{-\lambda_{1}s}\right)^{k}}\nonumber\\
&&\times\Big(\sum_{(u_{1},\cdots,u_{k})\subset Z_{s}}w_{-}\eta_{-}(A_{1})\frac{h(z_{u_{1}}(s))}{w(z_{u_{1}}(s))}\e^{-\lambda_{1}s}\times\cdots\times w_{-}\eta_{-}(A_{m})\frac{h(z_{u_{k}}(s))}{w(z_{u_{k}}(s))}\e^{-\lambda_{1}s}\Big)\nonumber\\
&\le&\frac{\prod_{i=1}^{m}\theta^{(i)}_{6}(t)^{k_{i}}}{\left(1-\theta_{5}(t)w_{-}\eta_{-}(A)c\e^{-\lambda_{1}s}\right)^{k}}\nonumber\\
&&\times\Big(\sum_{u_{1}\in Z_{s}}w_{-}\eta_{-}(A_{1})\frac{h(z_{u_{1}}(s))}{w(z_{u_{1}}(s))}\e^{-\lambda_{1}s}\Big)\times\cdots\times \Big(\sum_{u_{k}\in Z_{s}}w_{-}\eta_{-}(A_{m})\frac{h(z_{u_{k}}(s))}{w(z_{u_{k}}(s))}\e^{-\lambda_{1}s}\Big)\nonumber\\
&=&\frac{\prod_{i=1}^{m}\theta^{(i)}_{6}(t)^{k_{i}}}{\left(1-\theta_{5}(t)w_{-}\eta_{-}(A)c\e^{-\lambda_{1}s}\right)^{k}}\prod_{i=1}^{m}\left(w_{-}\eta_{-}(A_{i})W^{h/w}_{s}(Z)\right)^{k_{i}}.
\end{eqnarray}

For the lower bound, we have by \eqref{lem4.5.3} that on $\mathcal{E}^{1}_{t}$,
\begin{eqnarray}\label{prop4.12}
&&\sum_{(u_{1},\cdots,u_{k})\subset Z_{s}}\frac{\bp_{\delta_{z_{u_{1}}(s)}}\Big(Z_{t-s}(A_{1}+R(t))=1\Big)}{\bp_{\delta_{z_{u_{1}}(s)}}\Big(Z_{t-s}(A_{1}+R(t))=0\Big)}\times \cdots \times
\frac{\bp_{\delta_{z_{u_{k}}(s)}}\Big(Z_{t-s}(A_{m}+R(t))=1\Big)}{\bp_{\delta_{z_{u_{k}}(s)}}\Big(Z_{t-s}(A_{1}+R(t))=0\Big)}\nonumber\\
&\ge&\sum_{(u_{1},\cdots,u_{k})\subset Z_{s}}\bp_{\delta_{z_{u_{1}}(s)}}\Big(Z_{t-s}(A_{1}+R(t))=1\Big)\times \cdots \times
\bp_{\delta_{z_{u_{k}}(s)}}\Big(Z_{t-s}(A_{m}+R(t))=1\Big)\nonumber\\
&\ge&\sum_{(u_{1},\cdots,u_{k})\subset Z_{s}}\theta^{(1)}_{7}(t)w_{-}\eta_{-}(A_{1})\frac{h(z_{u_{1}}(s))}{w(z_{u_{1}}(s))}\e^{-\lambda_{1}s}
\times\cdots\times \theta^{(m)}_{7}(t)w_{-}\eta_{-}(A_{m})\frac{h(z_{u_{k}}(s))}{w(z_{u_{k}}(s))}\e^{-\lambda_{1}s}\nonumber\\
&=&\Big[\prod_{i=1}^{m}\left(\theta^{(i)}_{7}(t)w_{-}\eta_{-}(A_{i})\right)^{k_{i}}\Big]\times \Big[\sum_{(u_{1},\cdots,u_{k})\subset Z_{s}}\frac{h(z_{u_{1}}(s))}{w(z_{u_{1}}(s))}\e^{-\lambda_{1}s}
\times\cdots\times \frac{h(z_{u_{k}}(s))}{w(z_{u_{k}}(s))}\e^{-\lambda_{1}s}\Big].
\end{eqnarray}
Note that the sum
$\sum_{u_{1},\cdots,u_{k}\in Z_{s}}$ is no larger than the sum of $\sum_{1\le i<j\le k}\sum_{u_{1},\cdots,u_{k}\in Z_{s},u_{i}=u_{j}}$ and $\sum_{(u_{1},\cdots,u_{k})\subset Z_{s}}$,
and that
\begin{eqnarray}
&&\sum_{u_{1},\cdots,u_{k}\subset Z_{s},u_{i}=u_{j}}\frac{h(z_{u_{1}}(s))}{w(z_{u_{1}}(s))}\e^{-\lambda_{1}s}
\times\cdots\times \frac{h(z_{u_{k}}(s))}{w(z_{u_{k}}(s))}\e^{-\lambda_{1}s}\nonumber\\
&\le&c\e^{-\lambda_{1}s}\sum_{u_{1},\cdots,u_{k-1}\subset Z_{s}}\frac{h(z_{u_{1}}(s))}{w(z_{u_{1}}(s))}\e^{-\lambda_{1}s}
\times\cdots\times \frac{h(z_{u_{k}}(s))}{w(z_{u_{k}}(s))}\e^{-\lambda_{1}s}\nonumber\\
&\le&c\e^{-\lambda_{1}s}W^{h/w}_{s}(Z)^{k-1}.\nonumber
\end{eqnarray}
Thus we have
\begin{eqnarray}\nonumber
W^{h/w}_{s}(Z)^{k}&=&\sum_{u_{1},\cdots,u_{k}\subset Z_{s}}\frac{h(z_{u_{1}}(s))}{w(z_{u_{1}}(s))}\e^{-\lambda_{1}s}
\times\cdots\times \frac{h(z_{u_{k}}(s))}{w(z_{u_{k}}(s))}\e^{-\lambda_{1}s}\nonumber\\
&\le&\sum_{(u_{1},\cdots,u_{k})\subset Z_{s}}\frac{h(z_{u_{1}}(s))}{w(z_{u_{1}}(s))}\e^{-\lambda_{1}s}
\times\cdots\times \frac{h(z_{u_{k}}(s))}{w(z_{u_{k}}(s))}\e^{-\lambda_{1}s}\nonumber\\
&&\quad\quad+c\e^{-\lambda_{1}s}W^{h/w}_{s}(Z)^{k-1}.
\end{eqnarray}
Putting this back to \eqref{prop4.12}
we get that on $\mathcal{E}^{1}_{t}$,
\begin{eqnarray}\label{prop4.13}
&&\sum_{(u_{1},\cdots,u_{k})\subset Z_{s}}\frac{\bp_{\delta_{z_{u_{1}}(s)}}\Big(Z_{t-s}(A_{1}+R(t))=1\Big)}{\bp_{\delta_{z_{u_{1}}(s)}}\Big(Z_{t-s}(A_{1}+R(t))=0\Big)}\times \cdots \times
\frac{\bp_{\delta_{z_{u_{k}}(s)}}\Big(Z_{t-s}(A_{m}+R(t))=1\Big)}{\bp_{\delta_{z_{u_{k}}(s)}}\Big(Z_{t-s}(A_{1}+R(t))=0\Big)}\nonumber\\
&\ge&\Big[\prod_{i=1}^{m}\left(\theta^{(i)}_{7}(t)w_{-}\eta_{-}(A_{i})\right)^{k_{i}}\Big]\times\Big[W^{h/w}_{s}(Z)-c\e^{-\lambda_{1}s}W^{h/w}_{s}(Z)^{k-1}\Big].
\end{eqnarray}
We note that the final terms of \eqref{prop4.11} and \eqref{prop4.13} converges
to $\prod_{i=1}^{m}\left(w_{-}\eta_{-}(A_{i})W^{h/w}_{\infty}(Z)\right)^{k_{i}}$ almost surely.
Thus \eqref{prop4.7} follows immediately. Therefore we complete the proof.\qed

\bigskip

\begin{proposition}\label{prop5}
For every $\mu\in\mf$, $((Z_{t}\pm\sqrt{\lambda_{1}/2}\,t)_{t\ge 0},\p_{\mu})$ converges in distribution to a Cox process directed by
$w_{\pm}W^{h/w}_{\infty}(Z)\eta_{\pm}(dx)$,
where $W^{h/w}_{\infty}(Z)$ is the martingale limit of
$((W^{h/w}_{t}(Z))_{t\geq 0}, \p_{\mu})$.
\end{proposition}

\proof For any $f\in C^{+}_{c}(\R)$ and $\mu\in\mf$,
\begin{eqnarray}
\p_{\mu}\left[\e^{-\langle f,Z_{t}\pm\sqrt{\frac{\lambda_{1}}{2}}t\rangle}\right]
&=&\p_{\mu}\left[\p_{\mu}\left[\e^{-\langle f,Z_{t}\pm\sqrt{\frac{\lambda_{1}}{2}}t\rangle}|Z_{0}\right]\right]\nonumber\\
&=&\p_{\mu}\left[\prod_{u\in Z_{0}}\bp_{\delta_{z_{u}(0)}}\left[\e^{-\langle f,Z_{t}\pm\sqrt{\frac{\lambda_{1}}{2}}t\rangle}\right]\right]\nonumber\\
&=&\exp\left\{-\int_{\R}\left(1-\bp_{\delta_{x}}\left[\e^{-\langle f,Z_{t}\pm\sqrt{\frac{\lambda_{1}}{2}}t\rangle}\right]\right)w(x)\mu(dx)\right\}.
\end{eqnarray}
The final equality is because $(Z_{0},\p_{\mu})$ is a Poisson point process with intensity $w\mu$. Similarly one can prove that for every $\lambda\ge 0$,
$$\p_{\mu}\left[\e^{-\lambda W^{h/w}_{\infty}(Z)}\right]=\exp\left\{-\int_{\R}\left(1-\bp_{\delta_{x}}\left[\e^{-\lambda W^{h/w}_{\infty}(Z)}\right]\right)w(x)\mu(dx)\right\}.$$
Since by Proposition \ref{prop4} $\lim_{t\to+\infty}\bp_{\delta_{x}}\left[\e^{-\langle f,Z_{t}\pm\sqrt{\frac{\lambda_{1}}{2}}t\rangle}\right]=\bp_{\delta_{x}}\left[\e^{-w_{\pm} W^{h/w}_{\infty}(Z)\langle 1-\e^{-f},\eta_{\pm}\rangle}\right]$ for all $x\in\R$, we get by the bounded convergence theorem that
$$\lim_{t\to+\infty}\p_{\mu}\left[\e^{-\langle f,Z_{t}\pm\sqrt{\frac{\lambda_{1}}{2}}t\rangle}\right]=\p_{\mu}\left[\e^{- w_{\pm}W^{h/w}_{\infty}(Z)\langle 1-\e^{-f},\eta_{\pm}\rangle}\right].$$
Hence we prove this proposition.\qed

\bigskip

For $t\ge 0$, let $\max Z_{t}:=\max\{z_{u}(t):u\in Z_{t}\}$ be the maximum displacement of the skeleton branching diffusion.
\begin{proposition}\label{prop6}
For any $\mu\in\mathcal{M}_{c}(\R)$ and $y\in \R$,
\begin{eqnarray}\label{prop6.0}
&&\lim_{t\to+\infty}\p_{\mu}\left(\max Z_{t}-\sqrt{\frac{\lambda_{1}}{2}}t\le y\,\Big|\,\ W^{h}_{\infty}(X)>0\right)\nonumber\\
&=&\p_{\mu}\left[\exp\Big\{-\frac{w_{-}C_{-}}{\sqrt{2\lambda_{1}}}\e^{-\sqrt{2\lambda_{1}}y}W^{h}_{\infty}(X)\Big\}\,\Big|\,W^{h}_{\infty}(X)>0\right].
\end{eqnarray}
\end{proposition}
This implies that conditioned on $\{W^{h}_{\infty}(X)>0\}$, the maximal displacement of the skeleton branching diffusion centered by $\sqrt{\lambda_{1}/2}\,t$ converges in distribution to a randomly shifted Gumbel distribution.

\noindent\textbf{Proof of Proposition \ref{prop6}:}
Fix an arbitrary $y\in\R$.
If we set $A=A_{1}=(y,+\infty)$ and $k=k_{1}=0$ in \eqref{prop4.5}, then we get
$$
\bp_{\delta_{x}}\left(Z_{t}\Big(\sqrt{\frac{\lambda_{1}}{2}}t+y,+\infty\Big)=0\,\Big|\,\mathcal{G}_{s}\right)
\to \e^{-\frac{w_{-}C_{-}}{\sqrt{2\lambda_{1}}}\e^{-\sqrt{2\lambda_{1}}y}W^{h/w}_{\infty}(Z)}$$
in probability as $t\to+\infty$. It follows immediately that
$$\lim_{t\to+\infty}\bp_{\delta_{x}}\left(\max Z_{t}-\sqrt{\frac{\lambda_{1}}{2}}t\le y\right)=\bp_{\delta_{x}}\left[\exp\Big\{-\frac{w_{-}C_{-}}{\sqrt{2\lambda_{1}}}\e^{-\sqrt{2\lambda_{1}}y}W^{h/w}_{\infty}(Z)\Big\}\right].$$
Using this and the branching property,
we can show that
for any $\mu\in\mc$,
\begin{equation}\label{5.42}
\lim_{t\to+\infty}\p_{\mu}\left(\max Z_{t}-\sqrt{\frac{\lambda_{1}}{2}}t\le y\right)=\p_{\mu}\left[\exp\Big\{-\frac{w_{-}C_{-}}{\sqrt{2\lambda_{1}}}\e^{-\sqrt{2\lambda_{1}}y}W^{h}_{\infty}(X)\Big\}\right].
\end{equation}
We note that
\begin{eqnarray}
&&\p_{\mu}\left(\max Z_{t}-\sqrt{\frac{\lambda_{1}}{2}}t\le y,\ W^{h}_{\infty}(X)=0\right)\nonumber\\
&=&\p_{\mu}\left(\max Z_{t}-\sqrt{\frac{\lambda_{1}}{2}}t\le y,\ W^{h/w}_{\infty}(Z)=0\right)\nonumber\\
&=&\p_{\mu}\left(\|Z_{0}\|=0\right)+\p_{\mu}\left(\bp_{Z_{t}}\left(W^{h/w}_{\infty}(Z)=0\right);\|Z_{0}\|\not=0,\ \max Z_{t}-\sqrt{\frac{\lambda_{1}}{2}}t\le y\right).\label{5.39}
\end{eqnarray}
Noting \eqref{lem4.6.2}
one has $\bp_{Z_{t}}\left(W^{h/w}_{\infty}(Z)=0\right)=0$ $\p_{\mu}$-a.s. on $\{\|Z_{0}\|\not=0\}$. Thus the second term in the right hand side of \eqref{5.39}
equals $0$, and one gets
\begin{equation}\label{5.43}
\p_{\mu}\left(\max Z_{t}-\sqrt{\frac{\lambda_{1}}{2}}t\le y,\ W^{h}_{\infty}(X)=0\right)=\p_{\mu}\left(\|Z_{0}\|=0\right)=\e^{-\langle w,\mu\rangle}.
\end{equation}
Then we have
\begin{eqnarray}
&&\p_{\mu}\left(\max Z_{t}-\sqrt{\frac{\lambda_{1}}{2}}t\le y\,\Big|\,\ W^{h}_{\infty}(X)>0\right)-\p_{\mu}\left[\exp\Big\{-\frac{w_{-}C_{-}}{\sqrt{2\lambda_{1}}}\e^{-\sqrt{2\lambda_{1}}y}W^{h}_{\infty}(X)\Big\}\,\Big|\,W^{h}_{\infty}(X)>0\right]\nonumber\\
&=&\frac{\p_{\mu}\left(\max Z_{t}-\sqrt{\frac{\lambda_{1}}{2}}t\le y\right)-\e^{-\langle w,\mu\rangle}}{\p_{\mu}\left(W^{h}_{\infty}(X)>0\right)}
-\frac{\p_{\mu}\left[\exp\left\{-\frac{w_{-}C_{-}}{\sqrt{2\lambda_{1}}}\e^{-\sqrt{2\lambda_{1}}y}W^{h}_{\infty}(X)\right\}\right]-\p_{\mu}\left(W^{h}_{\infty}(X)=0\right)}
{\p_{\mu}\left(W^{h}_{\infty}(X)>0\right)}.\nonumber
\end{eqnarray}
Hence \eqref{prop6.0} follows by \eqref{5.42}.\qed

\begin{remark}\rm
One can order the positions of the particles alive at time $t$ in a non-increasing order: $R_{t,1}\ge R_{t,2}\ge\cdots\ge R_{t,\|Z_{t}\|}$. Then similarly as in Proposition \ref{prop6}, one can get the weak convergence of $(R_{t,1},R_{t,2},\cdots,R_{(t,n)})$.
\end{remark}

\subsection{Proofs of Theorem \ref{them2} and Theorem \ref{them3}}
The main idea of the proof for Theorem \ref{them2} is from \cite[Lemma 4.17]{Kallenberg}: Suppose $\xi_{1},\xi_{2},\cdots$ are Cox processes on $\R$ directed by some random measures $\eta_{1},\eta_{2},\cdots$. Then $\xi_{n}$ converges in distribution to some $\xi$ if and only if $\eta_{n}$ converges in distribution to some $\eta$, in which case $\xi$ is distributed as a Cox process directed by $\eta$.

\noindent\textbf{Proof of Theorem \ref{them2}:} Fix $\mu\in\mf$. In view of Proposition \ref{propskeleton}(iii), $(Z_{t}\pm\sqrt{\lambda_{1}/2}\,t,\p_{\mu})$ is distributed as a Cox process directed by $w(x\mp\sqrt{\lambda_{1}/2}\,t)\left(X_{t}\pm\sqrt{\lambda_{1}/2}\,t\right)(dx)$.
It then follows by \cite[Lemma 4.17]{Kallenberg} and Proposition \ref{prop5} that the latter converges in distribution to
$w_{\pm}W^{h}_{\infty}(X)\eta_{\pm}(dx)$.
This implies that for every $f\in C^{+}_{c}(\R)$, $\int_{\R}f(x)w(x\mp\sqrt{\lambda_{1}/2}\,t)\left(X_{t}\pm\sqrt{\lambda_{1}/2}\,t\right)(dx)$ converges in distribution to $w_{\pm}W^{h}_{\infty}(X)\langle f,\eta_{\pm}\rangle$.
Recall that for $x\ge M$, $w=w_{-}$ and for $x\leq -M$, $w=w_+$.
Note that for $t$ large enough such that
$x+\sqrt{\lambda_{1}/2}\,t\ge M$ and $x-\sqrt{\lambda_{1}/2}\,t\le -M$
for all $x\in \mbox{supp}f$,
$\int_{\R}f(x)w(x\mp\sqrt{\lambda_{1}/2}\,t)\left(X_{t}\pm\sqrt{\lambda_{1}/2}\,t\right)(dx)=w_{\pm}\langle f,X_{t}\pm\sqrt{\lambda_{1}/2}\,t\rangle$. Thus one gets that $\langle f,X_{t}\pm\sqrt{\lambda_{1}/2}\,t\rangle$ converges in distribution to $W^{h}_{\infty}(X)\langle f,\eta_{\pm}\rangle$. This implies that $X_{t}\pm\sqrt{\lambda_{1}/2}\,t$ converges in distribution to
$W^{h}_{\infty}(X)\eta_{\pm}(dx)$.\qed

\begin{remark}\label{rm:supremum}\rm
\begin{description}
\item{(i)} Theorem \ref{them2} implies that for any bounded and compactly supported measurable function $f$ on $\R$ whose set of discontinuous points has zero Lebesgue measure,
$\langle f,X_{t}\pm\sqrt{\lambda_{1}/2}\,t\rangle$ converges in distribution to $W^{h}_{\infty}(X)\langle f,\eta_{\pm}\rangle$. In particular for any compact set $B\subset\R$ whose boundary has zero Lebesgue measure, $X_{t}\left(\mp\sqrt{\lambda_{1}/2}\,t+B\right)$ converges in distribution to $W^{h}_{\infty}(X)\eta_{\pm}(B)$.
\item{(ii)}
We use $\max X_{t}$ to denote the supremum of the support of $X_{t}$, i.e., $\max X_{t}:=\sup\{x:X_{t}(x,+\infty)>0\}$.
Let $m>0$ and $y\in\R$. We have
\begin{equation}\label{5.40}
\pp_{\mu}\left(\max X_{t}-\sqrt{\frac{\lambda_{1}}{2}}t>y\right)\ge \pp_{\mu}\left(\Big\langle 1_{(y,y+m)},X_{t}-\sqrt{\frac{\lambda_{1}}{2}}t\Big\rangle >0\right).
\end{equation}
 Note that $\Big\langle 1_{(y,y+m)},X_{t}-\sqrt{\frac{\lambda_{1}}{2}}t\Big\rangle$ converges in distribution to $W^{h}_{\infty}(X)\eta_{-}(y,y+m)$.
Hence, letting $t\to+\infty$ in \eqref{5.40}, we get that
\begin{equation}\label{5.41}
\liminf_{t\to+\infty}\pp_{\mu}\left(\max X_{t}-\sqrt{\frac{\lambda_{1}}{2}}t>y\right)\ge\pp_{\mu}\left(W^{h}_{\infty}(X)>0\right).
\end{equation}
Note that
\begin{eqnarray*}
&&\pp_{\mu}\left(\max X_{t}-\sqrt{\frac{\lambda_{1}}{2}}t>y\,\Big|\,W^{h}_{\infty}(X)>0\right)\\
&=&\frac{\pp_{\mu}\left(\max X_{t}-\sqrt{\frac{\lambda_{1}}{2}}t>y\right)-\pp_{\mu}\left(\max X_{t}-\sqrt{\frac{\lambda_{1}}{2}}t>y,W^{h}_{\infty}(X)=0\right)}{\pp_{\mu}\left(W^{h}_{\infty}(X)>0\right)}\\
&=&\frac{\pp_{\mu}\left(\max X_{t}-\sqrt{\frac{\lambda_{1}}{2}}t>y\right)-\pp_{\mu}\left(\e^{-\langle w,X_{t}\rangle};\max X_{t}-\sqrt{\frac{\lambda_{1}}{2}}t>y\right)}{\pp_{\mu}\left(W^{h}_{\infty}(X)>0\right)}\\
&\ge&\frac{\pp_{\mu}\left(\max X_{t}-\sqrt{\frac{\lambda_{1}}{2}}t>y\right)-\pp_{\mu}\left(\e^{-\langle w,X_{t}\rangle}\right)}{\pp_{\mu}\left(W^{h}_{\infty}(X)>0\right)}\\
&=&\frac{\pp_{\mu}\left(\max X_{t}-\sqrt{\frac{\lambda_{1}}{2}}t>y\right)-\e^{-\langle w,\mu\rangle}}{\pp_{\mu}\left(W^{h}_{\infty}(X)>0\right)}.
\end{eqnarray*}
Hence by \eqref{5.41}, we have for any $y\in\R$,
$$\liminf_{t\to+\infty}\pp_{\mu}\left(\max X_{t}-\sqrt{\frac{\lambda_{1}}{2}}t>y\,\Big|\,W^{h}_{\infty}(X)>0\right)\ge 1-\frac{\e^{-\langle w,\mu\rangle}}{1-\e^{-\langle w,\mu\rangle}}>0.$$
So conditioned on $\{W^{h}_{\infty}(X)>0\}$, the distributions of $\{\max X_{t}-\sqrt{\lambda_{1}/2}\,t:t\ge 0\}$ are not tight.
This is very different from the behavior we observe in Proposition \ref{prop6} for the skeleton.
Loosely speaking, this is because the
range of the super-Brownian motion is much `larger' than that of
the embedded skeleton
\end{description}
\end{remark}

\bigskip

\noindent\textbf{Proof of Theorem \ref{them3}:} We take $\delta=\sqrt{\lambda_{1}/2}$ and $\mu\in\mc$. Suppose $\mbox{supp}\mu\subset [-k,k]$ for some $0<k<+\infty$.
We have
$$\pp_{\mu}\left[\mathcal{X}^{\delta t}_{t}\right]=\int_{\R}\Bp_{x}\left[e_{\beta}(t),|B_{t}|\ge \delta t\right]\mu(dx).$$
Since by Lemma \ref{lem:esti2} for $t$ large enough
$$\Bp_{x}\left[e_{\beta}(t),|B_{t}|\ge \delta t\right]\le \theta_{+}(t)\frac{1}{\sqrt{2\lambda_{1}}}(C_{+}+C_{-})h(x)\quad\forall x\in [-k,k],$$
where $\theta_{+}(t)\to 1$ as $t\to+\infty$, we have
$$\pp_{\mu}\left[\mathcal{X}^{\delta t}_{t}\right]\le\theta_{+}(t)\frac{1}{\sqrt{2\lambda_{1}}}(C_{+}+C_{-})\int_{\R}h(x)\mu(dx).$$
This implies that
$$\sup_{t\ge 0}\pp_{\mu}\left(\mathcal{X}^{\delta t}_{t}>\lambda\right)\le \sup_{t\ge 0}\frac{\pp_{\mu}\left[\mathcal{X}^{\delta t}_{t}\right]}{\lambda}\to 0\quad\mbox{ as }\lambda\to+\infty.$$
So the distributions of $\{\mathcal{X}^{\delta t}_{t}:t\ge 0\}$ are tight.

Applying similar argument as in the proof of Proposition \ref{prop4}, one can show that for any $x\in\R$, integers $m,n\ge 0$, integers $k_{1},\cdots,k_{m},l_{1},\cdots,l_{n}\ge 0$ and Borel sets $A_{1},\cdots,A_{m},B_{1},\cdots,B_{n}$ such that $\inf A_{i}>-\infty$ and $\sup B_{j}<+\infty$ for $i=1,\cdots,m$, $j=1,\cdots,n$,
\begin{eqnarray}
&&\bp_{\delta_{x}}\left(\bigcap_{i=1}^{m}\{Z_{t}(\delta t+A_{i})=k_{i}\},\bigcap_{j=1}^{n}\{Z_{t}(-\delta t+B_{j})=l_{j}\}\,|\,\mathcal{G}_{s}\right)\nonumber\\
&\to&
\exp\{-w_{-}W^{h/w}_{\infty}(Z)\sum_{i=1}^{m}\eta_{-}(A_{i})-w_{+}W^{h/w}_{\infty}(Z)\sum_{j=1}^{n}\eta_{+}(B_{j})\}\nonumber\\
&&\prod_{i=1}^{m}\frac{\left(w_{-}W^{h/w}_{\infty}(Z)\eta_{-}(A_{i})\right)
^{k_{i}}}{k_{i}!}\prod_{j=1}^{n}\frac{\left(w_{+}W^{h/w}_{\infty}(Z)\eta_{+}(B_{j})\right)
^{l_{j}}}{l_{j}!}\label{them3.1}
\end{eqnarray}
in probability as $t\to+\infty$. This implies that the point process
   $((Z_{t}-\delta t)+(Z_{t}+\delta t),\bp_{\delta_{x}})$
converges in distribution to a Cox process directed by
$W^{h/w}_{\infty}(Z)(w_{-}\eta_{-}(dx)+w_{+}\eta_{+}(dx))$.
Applying similar argument as in the proof of Theorem \ref{them2}, one can further show that
the random measure
$((X_{t}-\delta t)+(X_{t}+\delta t),\pp_{\mu})$
converges in distribution to
$W^{h}_{\infty}(X)(\eta_{-}(dx)+\eta_{+}(dx))$.
On the other hand,
by taking $n=m=1$ and $A_{1}=-B_{1}=[0,+\infty)$ in \eqref{them3.1}, one gets
\begin{eqnarray*}
&&\bp_{\delta x}\left(Z_{t}([\delta t,+\infty))=k_{1},\ Z_{t}((-\infty,-\delta t])=l_{1}\,|\,\mathcal{G}_{s}\right)\\
&\to&\e^{-\frac{1}{\sqrt{2\lambda_{1}}}W^{h/w}_{\infty}(Z)(w_{-}C_{-}+w_{+}C_{+})}\frac{\left(w_{-}W^{h/w}_{\infty}(Z)/\sqrt{2\lambda_{1}}\right)
^{k_{1}}}{k_{1}!}\frac{\left(w_{+}W^{h/w}_{\infty}(Z)/\sqrt{2\lambda_{1}}\right)
^{l_{1}}}{l_{1}!}
\end{eqnarray*}
in probability as $t\to+\infty$. Using similar computations as in the proof of Proposition \ref{prop5}, one gets that for all $\lambda_{1},\lambda_{2}\ge 0$,
\begin{eqnarray}\label{them3.2}
&&\lim_{t\to+\infty}\p_{\mu}\left[\e^{-\lambda_{1}Z_{t}([\delta t,+\infty))-\lambda_{2}Z_{t}((-\infty,-\delta t])}\right]\nonumber\\
&=&\p_{\mu}\left[\exp\left\{-\left(1-\e^{-\lambda_{1}}\right)\frac{1}{\sqrt{2\lambda_{1}}}W^{h}_{\infty}(X)w_{-}C_{-}-\left(1-\e^{-\lambda_{2}}\right)\frac{1}{\sqrt{2\lambda_{1}}}W^{h}_{\infty}(X)w_{+}C_{+}\right\}\right].
\end{eqnarray}
Recall that given $X_{t}$, $Z_{t}$ is a Poisson point process with intensity $wX_{t}$. Thus for $t$ sufficiently large,
\begin{eqnarray}
&&\p_{\mu}\left[\e^{-\lambda_{1}Z_{t}([\delta t,+\infty))-\lambda_{2}Z_{t}((-\infty,-\delta t])}\right]\nonumber\\
&=&\p_{\mu}\left[\exp\left\{-\Big\langle \left(1-\e^{-\lambda_{1}}\right)\1_{[\delta t,+\infty)}+\left(1-\e^{-\lambda_{2}}\right)\1_{(-\infty,-\delta t]},wX_{t}\Big\rangle\right\}\right]\nonumber\\
&=&\p_{\mu}\left[\exp\left\{-\left(1-\e^{-\lambda_{1}}\right)w_{-}X_{t}([\delta t,+\infty))-\left(1-\e^{-\lambda_{2}}\right)w_{+}X_{t}((-\infty,-\delta t])\right\}\right].\nonumber
\end{eqnarray}
Hence by \eqref{them3.2} we have
\begin{eqnarray}
&&\lim_{t\to+\infty}\p_{\mu}\left[\exp\left\{-\left(1-\e^{-\lambda_{1}}\right)w_{-}X_{t}([\delta t,+\infty))-\left(1-\e^{-\lambda_{2}}\right)w_{+}X_{t}((-\infty,-\delta t])\right\}\right]\nonumber\\
&=&\p_{\mu}\left[\exp\left\{-\left(1-\e^{-\lambda_{1}}\right)\frac{1}{\sqrt{2\lambda_{1}}}W^{h}_{\infty}(X)w_{-}C_{-}-\left(1-\e^{-\lambda_{2}}\right)\frac{1}{\sqrt{2\lambda_{1}}}W^{h}_{\infty}(X)w_{+}C_{+}\right\}\right].\nonumber
\end{eqnarray}
For $0\le \lambda< w_{-}\wedge w_{+}$, taking $\lambda_{1},\lambda_{2}$ such that $(1-\e^{-\lambda_{1}})w_{-}=(1-\e^{-\lambda_{2}})w_{+}=\lambda$, one gets that
\begin{equation}\label{them3.3}
\lim_{t\to+\infty}\pp_{\mu}\left[\e^{-\lambda \mathcal{X}^{\delta t}_{t}}\right]=\pp_{\mu}\left[\e^{-\lambda(C_{+}+C_{-})\frac{1}{\sqrt{2\lambda_{1}}}W^{h}_{\infty}(X)}\right].
\end{equation}

Suppose $(\mathcal{X}^{\delta t}_{t},\pp_{\mu})$ converges in distribution to $\xi$ along a subsequence $\{t_{n}:n\ge 1\}\subset [0,+\infty)$, for some random variable $\xi$.
Let $F_{1}$ and $F_{2}$ be the distribution functions of $\xi$ and $(C_{+}+C_{-})W^{h}_{\infty}(X)/\sqrt{2\lambda_{1}}$
respectively.
It suffices to show that $F_{1}=F_{2}$. Let $D_{1}$ be the set of continuous points of $F_{1}$. We note that $\{\mathcal{X}^{\delta t}_{t}\le x\}\subseteq \{X_{t}((\delta t,\delta t+N))+X_{t}((-\delta t -N,-\delta t))\le x\}$ for all $x,N\in\R$. Thus for any $y\in D_{1}$,
\begin{eqnarray*}
F_{1}(y)&=&\lim_{n\to+\infty}\pp_{\mu}\left(\mathcal{X}^{\delta t_{n}}_{t_{n}}\le y\right)\\
&\le& \limsup_{n\to+\infty}\pp_{\mu}\left(X_{t_{n}}((\delta t_{n},\delta t_{n}+N))+X_{t_{n}}((-\delta t_{n}-N,-\delta t_{n}))\le y\right)\\
&\le&\pp_{\mu}\left(W^{h}_{\infty}(X)\left(\eta_{-}((0,N))+\eta_{+}((-N,0))\right)\le y\right)\\
&=&F_{2}\left(\frac{\eta_{-}((0,+\infty))+\eta_{+}((-\infty,0))}{\eta_{-}((0,N))+\eta_{+}((-N,0))}\,y\right).
\end{eqnarray*}
By letting $N\to+\infty$, one gets that $F_{1}(y)\le F_{2}(y)$. If $F_{1}(y)<F_{2}(y)$ for some $y\in D_{1}$, then there is some $\epsilon>0$ such that $F_{1}(x)<F_{2}(x)$ for all $x\in (y,y+\epsilon)$. This yields that for any $\lambda>0$,
$$\mathrm{E}\left[\e^{-\lambda \xi}\right]-\pp_{\mu}\left[\e^{-\lambda(C_{+}+C_{-})\frac{1}{\sqrt{2\lambda_{1}}}W^{h}_{\infty}(X)}\right]=\lambda\int_{0}^{+\infty}\e^{-\lambda x}\left(F_{1}(x)-F_{2}(x)\right)dx<0,$$
which contradicts \eqref{them3.3}. Thus we have $F_{1}(x)=F_{2}(x)$ for all $x\in D_{1}$ and hence for all $x\in\R$.\qed

\appendix
\section{Appendix}

\begin{lemma}\label{lem4.1}
The martingale function $w$
 in \eqref{a2} is a solution to the following equation.
\begin{equation}\label{lemA0.1}
\frac{1}{2}w''(x)-\psi(x,w(x))=0,\quad\forall x\in \R.
\end{equation}
\end{lemma}

\proof It is proved by \cite[Lemma 2.1]{EKW} that the martingale function $w$ which satisfies \eqref{a2} is continuous on $\R$. Moreover, the argument leading to \cite[(2.4)]{EKW} shows that for any compact set $D$ of $\R$,
$$w(x)=\Bp_{x}\left[w(B_{t\wedge \tau_{D}})\right]-\Bp_{x}\left[\int_{0}^{t\wedge \tau_{D}}\psi(B_{s},w(B_{s}))ds\right],\quad\forall t\ge 0,\ x\in\R,$$
where $\tau_{D}$ denotes the first exit time of Brownian motion from $D$. Since $w$ is continuous and locally bounded, letting $t\to+\infty$ in the above equation, we get by the bounded convergence theorem that
$$w(x)=\Bp_{x}\left[w(B_{\tau_{D}})\right]-\Bp_{x}\left[\int_{0}^{\tau_{D}}\psi(B_{s},w(B_{s}))ds\right],\quad x\in D.$$
Applying similar argument as in the last paragraph of Page 708 in \cite{EP99}, one can show that $w$ is a solution to \eqref{lemA0.1}.\qed

\bigskip

\begin{lemma}\label{lemA1}
Suppose $\{\xi_{n}:n\ge 1\}$ is a sequence of point processes on $\R$, and $\eta$ is a locally finite random measure on $\R$. Then $\xi_{n}$ converges in distribution to a Cox process directed by $\eta$ if the following conditions hold.
\begin{description}
\item{(i)} For $m\in\mathbb{N}$, mutually disjoint bounded Borel sets $A_{1},\cdots,A_{m}$ of $\R$ and $k_{1},\cdots,k_{m}\in\mathbb{Z}^{+}$,
$$\mathrm{P}\left(\xi_{n}(A_{1})=k_{1},\cdots,\xi_{n}(A_{m})=k_{m}\right)\to \mathrm{E}\left[\e^{-\sum_{i=1}^{m}\eta(A_{i})}\prod_{i=1}^{m}\frac{\eta(A_{i})^{k_{i}}}{k_{i}!}\right]\quad\mbox{ as }n\to+\infty.$$
\item{(ii)} For any bounded Borel set $A$ of $\R$, $\sup_{n}\mathrm{E}\left[\xi_{n}(A)\right]<+\infty$.
\end{description}
\end{lemma}

\proof We need to show that for all $f\in C^{+}_{c}(\R)$,
\begin{equation}\label{lemA1.1}
\mathrm{E}\left[\e^{-\langle f,\xi_{n}\rangle}\right]\to \mathrm{E}\left[\e^{-\langle 1-\e^{-f},\eta\rangle}\right]\quad\mbox{ as }n\to+\infty.
\end{equation}
It is easy to deduce from (i) that \eqref{lemA1.1} holds if $f$ is a nonnegative compactly supported simple function. For an arbitrary $f\in C^{+}_{c}(\R)$ with supp$f\subset A$ where $A$ is a bounded Borel set of $\R$, one can find a nondecreasing sequence of nonnegative compactly supported simple functions $\{f_{k}:k\ge 1\}$ such that $f_{n}$ converges uniformly to $f$.
We note that for $k,n\ge 1$,
\begin{eqnarray*}
\left|\mathrm{E}\left[\e^{-\langle f_{k},\xi_{n}\rangle}\right]-\mathrm{E}\left[\e^{-\langle f,\xi_{n}\rangle}\right]\right|
&\le&\mathrm{E}\left[\left|\e^{-\langle f_{k},\xi_{n}\rangle}-\e^{-\langle f,\xi_{n}\rangle}\right|\right]\\
&\le&\mathrm{E}\left[\left|\langle f_{k},\xi_{n}\rangle-\langle f,\xi_{n}\rangle\right|\right]\\
&\le&\mathrm{E}\left[\langle |f_{k}-f|,\xi_{n}\rangle\right]\\
&\le&\|f_{k}-f\|_{\infty}\mathrm{E}\left[\xi_{n}(A)\right].
\end{eqnarray*}
It follows by (ii) that $\sup_{n}\left|\mathrm{E}\left[\e^{-\langle f_{k},\xi_{n}\rangle}\right]-\mathrm{E}\left[\e^{-\langle f,\xi_{n}\rangle}\right]\right|\to 0$ as $k\to+\infty$. So we have
\begin{eqnarray*}
\lim_{n\to+\infty}\mathrm{E}\left[\e^{-\langle f,\xi_{n}\rangle}\right]&=&\lim_{n\to+\infty}\lim_{k\to+\infty}\mathrm{E}\left[\e^{-\langle f_{k},\xi_{n}\rangle}\right]\\
&=&\lim_{k\to+\infty}\lim_{n\to+\infty}\mathrm{E}\left[\e^{-\langle f_{k},\xi_{n}\rangle}\right]\\
&=&\lim_{k\to +\infty}\mathrm{E}\left[\e^{-\langle 1-\e^{-f_{k}},\eta\rangle}\right]\\
&=&\mathrm{E}\left[\e^{-\langle 1-\e^{-f},\eta\rangle}\right].
\end{eqnarray*}
The first and final equalities are from the bounded convergence theorem.\qed

\medskip
{\bf Yan-Xia Ren}

LMAM School of Mathematical Sciences \& Center for
Statistical Science, Peking
University,

Beijing, 100871, P.R. China.

E-mail: yxren@math.pku.edu.cn

\medskip
{\bf Ting Yang}

School of Mathematics and Statistics, Beijing Institute of Technology,

Beijing, 100081, P.R.China;

Beijing Key Laboratory on MCAACI,

Beijing, 100081, P.R. China.

Email: yangt@bit.edu.cn

\end{document}